\documentclass[12pt,a4paper]{article} 
\usepackage{amsmath,amsfonts,amssymb,amsthm}
\usepackage{fullpage}
\title{Noncommutative instantons from \\ twisted conformal symmetries}
\author{Giovanni Landi$^1$, Walter D. van Suijlekom$^2$ \\[10mm]
$^1$ Dipartimento di Matematica e Informatica, Universit\`a di Trieste\\
Via A. Valerio 12/1, I-34127 Trieste, Italy \\
and INFN, Sezione di Trieste, Trieste, Italy\\
\texttt{landi@univ.trieste.it} \\[5mm] 
$^2$ Max Planck Institute for Mathematics\\
Vivatsgasse 7, D-53111 Bonn, Germany\\
\texttt{waltervs@mpim-bonn.mpg.de}}
\date{31 July 2006}

\renewcommand{\bar}[1]{\overline{#1}}
\renewcommand{\tilde}[1]{\widetilde{#1}}
\newcommand{\half}{{\mathchoice{\oh}{\oh}{\shalf}{\shalf}}} 
\newcommand{\oh}{{\tfrac{1}{2}}}    
\newcommand{\shalf}{{\scriptstyle\frac{1}{2}}} 
\def\II{\mathbb{I}}
\def\ii{\mathrm{i}} 
\def\into{\hookrightarrow}

\def\isom{\simeq}

\def\C{\mathbb{C}} 
\def\bH{\mathbb{H}}
\def\I{\mathbb{I}} 
 
\def\bR{\mathbb{R}}

\def\bT{\mathbb{T}}
\def\bZ{\mathbb{Z}} 
\def\A{\mathcal{A}}

\def\B{\mathcal{B}}
\def\cD{\mathcal{D}}
\def\E{\mathcal{E}}

\def\cH{\mathcal{H}}

\def\cO{\mathcal{O}}
 
\def\cS{\mathcal{S}} 
\def\cU{\mathcal{U}}

\def\cZ{\mathcal{Z}}
\def\g{\mathfrak{g}}

\def\SL{\mathrm{SL}}
\def\Sp{\mathrm{Sp}}

\def\chern{\mathrm{ch}}

\def\dix{\mathrm{Tr}_\omega}
\def\ind{\mathrm{index}}

\def\res{\mathrm{Res}}
\def\resz{\underset{z=0}{\res}}

\def\ad{\mathrm{ad}}

\def\Cinf{C^\infty}
\def\class{{(0)}}

\def\M{M_\theta}
\def\n{{(n)}}
\def\P{P_\theta}
\def\R{\bR_\theta}
\def\S{S_\theta}
\def\Sk{S_{\theta'}}

\def\T{\bT_\theta}

\def\Top{\mathrm{Top}}

\def\YM{\mathrm{YM}}
 
\def\lt{\triangleright}    
\def\rt{\triangleleft}    

\def\bd{\begin{displaymath} }
\def\ed{\end{displaymath} }
\def\be{\begin{equation}}
\def\ee{\end{equation}}
\def\bea{\begin{eqnarray}}
\def\eea{\end{eqnarray}}
\def\bmult{\begin{multline}}
\def\emult{\end{multline}}
\def\nn{\nonumber}
\newtheorem{thm}{Theorem}
\newtheorem{corl}[thm]{Corollary}
\newtheorem{lma}[thm]{Lemma} 
\newtheorem{prop}[thm]{Proposition}
\newtheorem{defn}[thm]{Definition}
\newtheorem{ex}[thm]{Example}
\newtheorem{rem}[thm]{\bf Remark}
\def\mb{\mbox{ }}

\newcommand{\bean}{\begin{eqnarray*}}
\newcommand{\eean}{\end{eqnarray*}}
\newcommand{\cinf}{C^\infty}       
\numberwithin{equation}{section}
\newcommand{\ca}{{\mathcal A}}

\newcommand{\ce}{{\mathcal E}}

\newcommand{\ch}{{\mathcal H}}

\newcommand{\cs}{{\mathcal S}}

\newcommand{\cu}{{\mathcal U}}

\newcommand{\IC}{{\mathbb C}}

\newcommand{\IR}{{\mathbb R}}

\newcommand{\IT}{{\mathbb T}}
\newcommand{\IZ}{{\mathbb Z}}

\def\lra{\longrightarrow}

\newcommand{\wh}{\widehat}
\newcommand{\wt}{\widetilde}
\def\bar#1{\overline{#1}}

\newcommand{\omca}{\Omega {\mathcal A}}

\newcommand{\oca}[1]{\Omega^{#1}{\mathcal A}}

\newcommand{\comca}{{\mathcal E}\otimes_{\mathcal A}\Omega {\mathcal A}}
\newcommand{\coca}[1]{{\mathcal E}\otimes_{\mathcal A}\Omega^{#1}{\mathcal A}}

\newcommand{\ota}{\otimes_{\mathcal A}}

\newcommand{\otc}{\otimes_{\IC}}
\newcommand{\ot}{\otimes}

\newcommand{\op}{\oplus}
\newcommand{\hs}[2]{\left\langle #1,#2\right\rangle}

\DeclareMathOperator{\id}{id}

\DeclareMathOperator{\Tr}{Tr}

\DeclareMathOperator{\tr}{tr}
\DeclareMathOperator{\SU}{SU}

\DeclareMathOperator{\U}{U} 
\DeclareMathOperator{\SO}{SO}
\DeclareMathOperator{\Spin}{Spin}
\DeclareMathOperator{\Aut}{Aut} 
\DeclareMathOperator{\Hom}{Hom}
\DeclareMathOperator{\End}{End}
\newcommand{\dd}{{\rm d}} 
\newbox\ncintdbox \newbox\ncinttbox
\setbox0=\hbox{$-$} \setbox2=\hbox{$\displaystyle\int$}
\setbox\ncintdbox=\hbox{\rlap{\hbox to
\wd2{\hskip-.125em\box2\relax \hfil}}\box0\kern.1em}
\setbox0=\hbox{$\vcenter{\hrule width 4pt}$}
\setbox2=\hbox{$\textstyle\int$}
\setbox\ncinttbox=\hbox{\rlap{\hbox to
\wd2{\hskip-.175em\box2\relax \hfil}}\box0\kern.1em}
\newcommand{\ncint}{\mathop{\mathchoice{\copy\ncintdbox}%
{\copy\ncinttbox}{\copy\ncinttbox}{\copy\ncinttbox}}\nolimits}

\begin{document}

\maketitle
\hyphenation{de-for-ma-tions}
\begin{abstract}
 We construct a five-parameter family of gauge-nonequivalent $SU(2)$ instantons 
 on a noncommutative four sphere $\S^4$ and of topological charge equal to $1$. 
 These instantons are critical points of a gauge functional and satisfy 
 self-duality equations with respect to a Hodge star operator on forms on $\S^4$. 
 They are obtained by acting with a twisted conformal symmetry on a basic 
 instanton canonically associated  with a  noncommutative instanton bundle on the 
 sphere. A completeness argument for this family is obtained by means of index 
 theorems. The dimension of the ``tangent space'' to the moduli space  is 
 computed as the index of a twisted Dirac operator and turns out to be equal to 
 five,  a number that survives deformation. 
\end{abstract}

\thispagestyle{empty}
\newpage
\tableofcontents
\newpage

\section{Introduction}
 The importance of Yang--Mills instantons in physics and mathematics needs not be 
 stressed. They have played a central role since their first appearance 
 \cite{BPST} and are most elegantly described via the so-called ADHM construction 
 \cite{ADHM78,Ati79}. The generalization in \cite{NS98} of this method for 
 instantons on a noncommutative space $\IR^4$ has found several important 
 applications notably in brane and superconformal theories.

 Toric noncommutative manifolds $M_\Theta$ were constructed and studied in 
 \cite{CL01}. One starts with any (Riemannian spin) manifold $M$ carrying a torus 
 action and then deforms the torus to a noncommutative one governed by a real 
 antisymmetric matrix $\Theta$ of deformation parameters. 
 The starting example of \cite{CL01} -- the archetype of all these deformations 
 -- was a four dimensional sphere $\S^4$,
 which came with a natural noncommutative instanton bundle endowed with a natural 
 connection. At the classical value of the deformation parameter, $\theta=0$, the 
 bundle and the connection reduces to the one of \cite{BPST}. 
 The present sphere $\S^4$ can be thought of \cite{CD02} as a one point 
 compactification of a noncommutative $\R^4$ which is structurally different from 
 the one considered in \cite{NS98}.

In \cite{LS04} this basic noncommutative instanton was put in the context of an 
$\SU(2)$ noncommutative principal fibration $\Sk^7 \to \S^4$ over $\S^4$ . In the 
present paper, we continue the analysis and consider it in the setting of a 
noncommutative Yang--Mills theory. We then construct a five-parameter family of 
(infinitesimal) gauge-nonequivalent instantons, by acting with twisted conformal 
symmetries on the basic instanton. All these instantons will be gauge 
configurations satisfying self-duality equations -- with a suitable defined Hodge 
$\ast_\theta$-operator on forms $\Omega(\S^4)$ -- and will have a ``topological 
charge'' of value $1$. A completeness argument on the family of instantons is 
provided by index theoretical arguments, similar to the one in \cite{AHS78} for 
undeformed instantons on $S^4$. The dimension of the ``tangent'' of the moduli 
space can be computed as the index of a twisted Dirac operator which turns out to 
be equal to its classical value that is five. The twisting of the conformal 
symmetry is implemented with a twist of Drinfel'd type \cite{Dr83,Dr90} -- in 
fact, explicitly constructed by Reshetikhin \cite{Re90} -- and gives rise to a 
deformed Hopf algebra $\U_\theta(so(5,1))$. That these are conformal 
infinitesimal transformations is stressed by the fact that the Hopf algebra 
$\U_\theta(so(5,1))$ leaves the Hodge $\ast_\theta$-structure of $\Omega(\S^4)$ 
invariant.

The paper is organized as follows. In Section~\ref{section:connections} we recall 
the setting of gauge theories (connections) and gauge transformations on finite 
projective modules (the substitute for vector bundles) over algebras (the 
substitute for spaces). The main objective is to implement a Bianchi identity 
that will be crucial later on for the self-duality equations.

Section~\ref{section:toric-ncm} deals with toric noncommutative manifolds. These 
were indeed named isospectral deformations in that they can be endowed with the 
structure of a noncommutative Riemannian spin manifold via a spectral triple 
$(C^\infty(M_\theta), D, \ch$) with the properties of \cite{C96}. For this class 
of examples, the Dirac operator $D$ is the classical one and $\ch=L^2(M,\cs)$ is 
the usual Hilbert space of spinors
 on which the algebra $C^\infty(M_\theta)$ acts in a twisted manner. Thus one 
twists the algebra and its representation while keeping the geometry unchanged.  
The resulting noncommutative geometry is
 isospectral and all spectral properties are preserved including the dimension. 
Both the algebra and its action on spinors can be given via a ``star-type" 
product.

In Section~\ref{section:gts} we specialize to gauge theories on the sphere $\S^4$ 
and introduce a Yang--Mills action functional, from which we derive field 
equations (equations for critical points), as well as a topological action 
functional whose absolute value gives a lower bound for the Yang--Mills action.

The heart of the paper is Section~\ref{sect:constrinst} were we explicitly 
construct instantons. As usual, these are gauge configurations which are 
solutions of (anti)self-duality equations and realize absolute minima of the 
Yang--Mills functional. We start from a basic instanton which is shown to be 
invariant under twisted orthogonal transformations in $\U_\theta(so(5))$. We then 
perturb it by the action of conformal operators in $\U_\theta(so(5,1)) -
 \U_\theta(so(5))$ producing a five parameter family of new, not gauge equivalent 
instantons. A completeness argument is obtained by using an index theorem to 
compute the dimension of the tangent space of the moduli space of instantons on 
$\S^4$, which is shown to be just five. The relevant material from noncommutative 
index theory is recalled in the appendix.
 
Section~\ref{section:general} sketches a general scheme for gauge theories on 
four dimensional toric noncommutative manifolds.

\section{Connections and gauge transformations}\label{section:connections}

We first review the notion of a (gauge) connection on a (finite projective) 
module $\ce$ over an algebra $\ca$ with respect to a given calculus; we take a 
right module structure. Also, we recall gauge transformations in this setting. We 
refer to \cite{C94} for more details (see also \cite{Lnd97}).

\subsection{Connections on modules}

Let us suppose we have an algebra $\ca$ with a differential calculus 
$(\omca=\op_p \oca{p}, \dd)$. A {\it connection} on the right $\ca$-module $\ce$ 
is a $\IC$-linear map \[ \nabla : \coca{p} \lra \coca{p+1} , \] defined for any 
$p \geq 0$, and satisfying the Leibniz rule \[ \nabla(\omega \rho) = (\nabla 
\omega) \rho + (-1)^{p} \omega \dd \rho , ~~\forall ~\omega \in \coca{p} , ~\rho 
\in \omca .  \] \noindent A connection is completely determined by its 
restriction \be \nabla : \ce \to \coca{1} , \ee
 which satisfies \be\label{ulei} \nabla (\eta a) = (\nabla \eta) a + \eta \ota 
\dd a , ~~\forall ~\eta \in \ce , ~a \in \ca , \ee and which is extended to all 
of $\coca{p}$ using Leibniz rule. It is the latter rule that implies the 
$\omca$-linearity of the composition, \[ \nabla^2 = \nabla \circ \nabla : 
\coca{p} \lra \coca{p+2} . \] Indeed, for any $\omega \in \coca{p}$ and $\rho \in 
\omca$ it follows that $ \nabla^2 ({\omega \rho}) = \nabla \big( (\nabla \omega) 
\rho + (-1)^{p} \omega \dd \rho \big)  = (\nabla ^2 \omega) \rho + (-1)^{p+1} 
(\nabla \omega) \dd \rho + (-1)^{p} (\nabla \omega) \dd \rho + \omega \dd^2 \rho 
= (\nabla ^2 \omega) \rho $. The restriction of $\nabla^2$ to $\ce$ is the {\it 
curvature} \be F : \ce \to \coca{2} , \ee of the connection. It is $\ca$-linear, 
$F(\eta a) = F(\eta)a$ for any $\eta\in\ce, a \in \ca$, and satisfies \be 
\nabla^2(\eta \ota \rho) = F(\eta) \rho , ~~\forall ~\eta \in \ce , ~\rho \in 
\omca . \ee Thus, $ F\in \Hom_{\ca}(\ce, \coca{2} )$, the latter being the 
collection of (right) $\A$-linear homomorphisms from $\E$ to $\coca{2}$ (an 
alternative notation for this collection that is used in the literature, is 
$\End_{\ca}(\ce, \coca{2} )$).

In order to have the notion of a Bianchi identity we need some generalization. 
Let $\End_{\omca}(\comca)$ be the collection of all $\omca$-linear endomorphisms 
of $\comca$. It is an algebra under composition. The curvature $ F$ can be 
thought of as an element of $\End_{\omca}(\comca)$. There is a map \begin{align} 
\label{def:conn-end} [\nabla, ~\cdot~ ] & ~:~ \End_{\omca}(\comca) \lra 
\End_{\omca}(\comca) , \nn \\ [\nabla, T] &:= \nabla \circ T - (-1)^{|T|}~ T 
\circ \nabla , \end{align} where $|T|$ denotes the degree of $T$ with respect to 
the $\bZ^2$-grading of $\Omega\A$. Indeed, for any $\omega \in \coca{p}$ and 
$\rho \in \omca$, it follows that \begin{align*} [\nabla, T] ({\omega \rho}) & = 
\nabla (T ( \omega \rho) ) - (-1)^{|T|} ~ T ( \nabla ( \omega \rho)) \nn \\ & = 
\nabla \big( T ( \omega) \rho\big) - (-1)^{|T|} ~ T \big( (\nabla \omega) \rho + 
(-1)^{p} \omega \dd \rho \big) \nn \\ & = \big(\nabla (T(\omega))\big) \rho + 
(-1)^{p+|T|}~ T(\omega) \dd \rho - (-1)^{|T|}~ T (\nabla \omega)  \rho - 
(-1)^{p+|T|}~ T(\omega) \dd \rho \nn \\ & = \big(\nabla (T(\omega)) - (-1)^{|T|} 
~T (\nabla \omega) \big) \rho = \big([\nabla, T] ( \omega)  \big) \rho , 
\end{align*} and the map in \eqref{def:conn-end} is well-defined. It is 
straightforwardly checked that $[\nabla, ~\cdot~]$ is a graded derivation for the 
algebra $\End_{\omca}(\comca)$, \be [\nabla,S \circ T]=[\nabla,S]\circ T + 
(-1)^{|S|} S \circ [\nabla,T]. \ee \begin{prop}\label{ubianchi} The curvature $ 
F$ satisfies the {\rm Bianchi identity}, \be [\nabla, F ] = 0 .  \label{ubia} \ee 
\end{prop} \begin{proof} Since $ F$ is an even element in $\End_{\omca}(\comca)$, 
the map $[\nabla, F ]$ makes sense. Furthermore, \[ [\nabla, F ] = \nabla \circ 
\nabla^2 - \nabla^2 \circ \nabla = \nabla^3 - \nabla^3 = 0 . \] \end{proof} 
\noindent In Section~II.2 of \cite{C85}, such a Bianchi identity was implicitly 
used in the construction of a so-called canonical cycle from a connection on a 
finite projective $\A$-module $\E$.

Connections always exist on a projective module. On the module $\ce = \IC^N \otc 
\ca \simeq \ca^N$, which is free, a connection is given by the operator \[ 
\nabla_0 = \II \ot \dd : \IC^N \otc \oca{p} \lra \IC^N \otc \oca{p+1} . \] With 
the canonical identification $\IC^N \otc \omca = (\IC^N \otc \ca) \ota \omca 
~\simeq~ (\omca)^N$, one thinks of $\nabla_0$ as acting on $(\omca)^N$ as the 
operator $\nabla_0 = (\dd, \dd, \cdots, \dd)$ ($N$-times).  Next, take a 
projective module $\ce$ with inclusion map, $\lambda : \ce \to \ca^N$, which 
identifies $\ce$ as a direct summand of the free module $\ca^N$ and idempotent $p 
: \ca^N \to \ce$ which allows one to identify $\ce = p \ca^N$. Using these maps 
and their natural extensions to $\ce$-valued forms, a connection $\nabla_0$ on 
$\ce$ (called {\it Levi-Civita} or {\it Grassmann}) is the composition, \[ 
\coca{p} ~\stackrel{\lambda}{\lra}~ \IC^N \otc \oca{p} ~\stackrel{\II \ot 
\dd}{\lra}~ \IC^N \otc \oca{p+1} ~\stackrel{p}{\lra}~ \coca{p+1}, \] that is 
\be\label{ugras} \nabla_0 = p \circ (\II \ot \dd ) \circ \lambda.  \ee One 
indicates it simply by $\nabla_0 = p \dd$.  The space $C(\ce)$ of all connections 
on $\ce$ is an affine space modeled on $\Hom_{\ca}({\ce,\coca{1}})$. Indeed, if 
$\nabla_1, \nabla_2$ are two connections on $\ce$, their difference is 
$\ca$-linear, \[ (\nabla_1 - \nabla_2)(\eta a) = ((\nabla_1 - \nabla_2)(\eta )) a 
, \quad \forall ~ \eta \in \ce , ~a \in \ca , \] so that $\nabla_1 - \nabla_2 \in 
\Hom_{\ca}({\ce,\coca{1}})$. Thus, any connection can be written as \be \nabla = 
p \dd + \alpha , \label{uconn} \ee where $\alpha$ is any element in 
$\Hom_{\ca}({\ce,\coca{1}})$.  The ``matrix of $1$-forms'' $\alpha$ as in 
\eqref{uconn} is called the {\it gauge potential} of the connection $\nabla$. The 
corresponding curvature $ F$ of $\nabla$ is \be
 F = p \dd p \dd p + p \dd \alpha + \alpha^2 . \label{ucurv} \ee

\bigskip Next, let the algebra $\ca$ have an involution $^*$; it is extended to 
the whole of $\omca$ by the requirement $(\dd a)^* = \dd a^*$ for any $a\in\ca$. 
A {\it Hermitian structure} on the module $\ce$ is a map $\hs{ \cdot}{\cdot} : 
\ce \times \ce \to \ca$ with the properties \begin{align} \label{hest} &\hs{\eta 
}{\xi a } = \hs{\xi}{\eta} a , \qquad \hs{\eta}{\xi }^* = \hs{\xi}{\eta} , \nn \\ 
&\hs{\eta}{\eta} \geq 0 , \qquad \hs{\eta}{\eta} = 0 \iff \eta = 0 , \end{align} 
for any $\eta,\xi\in\ce$ and $a\in\ca$ (an element $a\in\ca$ is positive if it is 
of the form $a=b^* b$ for some $b\in\ca$). We shall also require the Hermitian 
structure to be {\it self-dual}, {i.e.} every right $\ca$-module homomorphism 
$\phi: \E \to \ca$ is represented by an element of $\eta \in \E$, by the 
assignment $\phi(\cdot) = \hs{\eta}{\cdot}$, the latter having the correct 
properties by the first of \eqref{hest}.

The Hermitian structure is naturally extended to an $\omca$-valued linear map on 
the product $\comca \times \comca$ by \be\label{uhext} \hs{ \eta \ota \omega}{\xi 
\ota \rho }~ = (-1)^{|\eta| |\omega|} \omega^* \hs{\eta}{\xi}\rho , \quad \forall 
~\eta, \xi \in \comca , ~\omega,\rho\in \omca. \ee

A connection $\nabla$ on $\ce$ and a Hermitian structure $\hs{ \cdot}{\cdot}$ on 
$\ce$ are said to be compatible if the following condition is satisfied 
\cite{C94}, \be\label{ucompat} \hs{\nabla \eta}{\xi} + \hs{\eta}{\nabla \xi} = 
\dd \hs{\eta}{\xi}, \quad \forall ~\eta, \xi \in \ce . \ee It follows directly 
from the Leibniz rule and \eqref{uhext} that this extends to 
\be\label{ucompat-ext} \hs{\nabla \eta}{\xi} + (-1)^{|\eta|} \hs{\eta}{\nabla 
\xi} = \dd \hs{\eta}{\xi} , \quad \forall ~\eta, \xi \in \comca . \ee On the free 
module $\A^N$ there is a canonical Hermitian structure given by 
\be\label{hs-canonical} \hs{\eta}{\xi} = \sum_{j=1}^N \eta_j^* \xi_j , \ee with 
$\eta = (\eta_1, \cdots, \eta_N)$ and $\eta = (\eta_1, \cdots, \eta_N)$ any two 
elements of $\A^N$.

Under suitable regularity conditions on the algebra $\A$ all Hermitian structures 
on a given finite projective module $\E$ over $\A$ are isomorphic to each other 
and are obtained from the canonical structure \eqref{hs-canonical} on $\A^N$ by 
restriction \cite[II.1]{C94}. Moreover, if $\E=p \A^N$, then $p$ is self-adjoint: 
$p=p^*$, with $p^*$ obtained by the composition of the involution $^*$ in the 
algebra $\A$ with the usual matrix transposition. The Grassmann connection 
\eqref{ugras} is easily seen to be compatible with this Hermitian structure, \be 
\dd \hs{\eta}{\xi} = \hs{\nabla_0 \eta}{\xi} + \hs{\eta}{\nabla_0 \xi} . \ee For 
a general connection \eqref{uconn}, the compatibility with the Hermitian 
structure reduces to \be \hs{ \alpha \eta}{\xi } + \hs{\eta}{\alpha \xi }= 0 , 
\quad \forall ~\eta, \xi \in \ce , \ee which just says that the gauge potential 
is skew-hermitian, \be \alpha ^* = -\alpha . \ee We still use the symbol $C(\ce)$ 
to denote the space of compatible connections on $\ce$.

Let $\End^s_{\omca}(\comca)$ denote the space of elements $T$ in 
$\End_{\omca}(\comca)$ which are skew-hermitian with respect to the Hermitian 
structure \eqref{uhext}, {i.e.} satisfying \be \label{def:skew} \hs{T \eta}{\xi} 
+ \hs{\eta}{T\xi}=0, \qquad \forall ~ \eta,\xi\in \E. \ee \begin{prop} 
\label{prop:conn-skew} The map $[\nabla,~\cdot~]$ in \eqref{def:conn-end} 
restricts to $\End^s_{\omca}(\comca)$ as a derivation \begin{align} [\nabla, 
~\cdot~ ] & ~:~ \End^s_{\omca}(\comca) \lra \End^s_{\omca}(\comca) , \end{align} 
\end{prop} \begin{proof} Let $T\in\End^s_{\omca}(\comca)$ be of order $|T|$; it 
then satisfies \be\label{inter} \hs{T \eta}{\xi} + (-1)^{|\eta||T|} 
\hs{\eta}{T\xi}=0, \ee for $\eta,\xi\in \comca$. Since $[\nabla,T]$ is 
$\omca$-linear, it is enough to show that $$ \hs{[\nabla,T]\eta}{\xi}+ 
\hs{\eta}{[\nabla,T]\xi}=0, \qquad \forall ~\eta,\xi \in \E. $$ This follows from 
equations \eqref{inter} and \eqref{ucompat-ext}, \begin{align*} 
\hs{[\nabla,T]\eta}{\xi}+ \hs{\eta}{[\nabla,T]\xi} &= \hs{\nabla T \eta}{\xi}- 
(-1)^{|T|} \hs{T \nabla \eta}{\xi}+\hs{\eta}{\nabla T \xi} - (-1)^{|T|} 
\hs{\eta}{T\nabla \xi} \\ &= \hs{\nabla T \eta}{\xi}- \hs{\nabla \eta}{T 
\xi}+\hs{\eta}{\nabla T \xi} - (-1)^{|T|} \hs{T\eta}{\nabla \xi} \\ &=\dd \big( 
\hs{T\eta}{\xi}+\hs{\eta}{T\xi} \big) = 0. \end{align*} \end{proof}

\subsection{Gauge transformations} We now add the additional requirement that the 
algebra $\A$ is a Fr\'echet algebra and that $\E$ a right Fr\'echet module. That 
is, both $\A$ and $\E$ are complete in the topology defined by a family of 
seminorms $\| \cdot \|_i$ such that the following condition is satisfied: for all 
$j$ there exists a constant $c_j$ and an index $k$ such that \be \| \eta a \|_j 
\leq c_j \| \eta \|_k \|a\|_k. \ee The collection $\End_{\ca}(\ce)$ of all 
$\ca$-linear maps is an algebra with involution; its elements are also called 
endomorphisms of $\ce$. It becomes a Fr\'echet algebra with the following family 
of seminorms: for $T \in \End_{\ca}(\ce)$, \be \| T \|_i = \sup_\eta \left\{ \| T 
\eta \|_i : \| \eta\|_i \leq 1 \right\}. \ee Since we are taking a self-dual 
Hermitian structure (see the discussion after \eqref{hest}), any $T \in 
\End_{\ca}(\ce)$ is adjointable, that is it admits an adjoint, an $\ca$-linear 
map $T^*:\ce\to\ce$ such that\[ \hs{T^* \eta}{\xi} = \hs{\eta}{T \xi} , \quad 
\forall ~\eta, \xi \in \ce . \] The group $\cu(\ce)$ of unitary endomorphisms of 
$\ce$ is given by \be\label{unit} \cu(\ce) := \{ u \in \End_{\ca}(\ce) ~|~ u u^* 
= u^* u = \id_\ce \} . \ee This group plays the role of the {\it infinite 
dimensional group of gauge transformations}. It naturally acts on compatible 
connections by \be\label{ugcon} (u, \nabla) \mapsto \nabla^u := u^* \nabla u , 
\quad \forall ~u\in \cu(\ce), ~\nabla \in C(\ce) , \ee where $u^*$ is really $u^* 
\ot \id_{\omca}$; this will always be understood in the following. Then the 
curvature transforms in a covariant way \be (u, F) \mapsto F^u = u^* F u , 
\label{ugcur} \ee since, evidently, $ F^u = (\nabla^u)^2 = u^* \nabla u u^* 
\nabla u^* = u^* \nabla^2 u = u^* F u$. \\ As for the gauge potential, one has 
the usual affine transformation, \be (u, \alpha) \mapsto \alpha^u := u^* p \dd u 
+ u^* \alpha u . \label{ugpot} \ee Indeed, $\nabla^u (\eta) = u^*( p \dd + 
\alpha) u \eta = u^* p \dd (u \eta) + u^* \alpha u \eta = u^* p u \dd \eta + u^* 
p (\dd u) \eta + u^* \alpha u \eta = p \dd \eta + (u^* p \dd u + u^* \alpha u) 
\eta$ for any $\eta\in\ce$, which yields \eqref{ugpot} for the transformed 
potential.

\bigskip The ``tangent vectors'' to the gauge group $\cU(\E)$ constitute the 
vector space of infinitesimal gauge transformations. Suppose $\{u_t\}_{t \in 
\bR}$ is a differentiable family of elements in $\End_\A(\E)$ (in the topology 
defined by the above sup-norms) and define $X:=(\partial u_t /\partial t)_{t=0}$. 
Unitarity of $u_t$ then induces that $X=-X^*$. In other words, for $u_t$ to be a 
gauge transformation, $X$ should be a skew-hermitian endomorphisms of $\E$. In 
this way, we understand $\End^s_\A(\E)$ as the collection of {\it infinitesimal 
gauge transformations}.  It is a real vector space whose complexification 
$\End^s_\A(\E)\ot_\bR \C$ can be identified with $\End_\A(\E)$.

Infinitesimal gauge transformations act on a connection in a natural way. Let the 
gauge transformation $u_t$, with $X=(\partial u_t /\partial t)_{t=0}$, act on 
$\nabla$ as in \eqref{ugcon}. From the fact that $(\partial (u_t \nabla u_t^* 
)/\partial t)_{t=0} = [\nabla, X]$, we conclude that an element $X \in 
\End^s_\A(\E)$ acts infinitesimally on a connection $\nabla$ by the addition of 
$[\nabla,X]$, \be\label{uXcon} (X, \nabla) \mapsto \nabla^X = \nabla + t 
[\nabla,X] + \cO(t^2) , \quad \forall ~X \in \End^s_\A(\E), ~\nabla \in C(\ce). 
\ee As a consequence, for the transformed curvature one finds \be\label{uXcur} 
(X, F) \mapsto F^X = F + t [F,X] + \cO(t^2), \ee since $F^X = (\nabla + t 
[\nabla,X])\circ (\nabla + t [\nabla,X]) = \nabla^2 + t [\nabla^2, X] + 
\cO(t^2)$.

\section{Toric noncommutative manifolds $\M$}\label{section:toric-ncm}

We start by recalling the general construction of toric noncommutative manifolds 
given in \cite{CL01} where they were called isospectral deformations. These are 
deformations of a classical Riemannian manifold and satisfy all the properties of 
noncommutative spin geometry~\cite{C96}. They are the content of the following 
result taken from~\cite{CL01},

\begin{thm}\label{Theorem6} Let $M$ be a compact spin Riemannian manifold whose 
isometry group has rank $r \geq 2$. Then $M$ admits a natural one parameter 
isospectral deformation to noncommutative geometries $M_{\theta}$. \end{thm}

The idea of the construction is to deform the standard spectral triple describing 
the Riemannian geometry of $M$ along a torus embedded in the isometry group, thus 
obtaining a family of spectral triples describing noncommutative geometries.

\subsection{Deforming a torus action} \label{subsection:def-torus} Let $M$ be an 
$m$ dimensional compact Riemannian manifold equipped with an isometric smooth 
action $\sigma$ of an $n$-torus $\bT^n$, $n \geq 2$.  We denote by $\sigma$ also 
the corresponding action of $\bT^n$ by automorphisms -- obtained by pull-backs -- 
on the algebra $\Cinf(M)$ of smooth functions on $M$. 

The algebra $\cinf(M)$ may be decomposed into spectral subspaces which are 
indexed by the dual group $\IZ^n = \wh\IT^n$. Now, with $s=(s_1, \cdots, 
s_n)\in\bT^n$, each $r \in \IZ^n$ yields a character of $\IT^n$, $e^{2\pi i s} 
\mapsto e^{2\pi i r\cdot s}$, with the scalar product $r\cdot s := r_1 s_1 
+\cdots+ r_n s_n$. The $r$-th spectral subspace for the action $\sigma$ of 
$\IT^n$ on $\cinf(M)$ consists of those smooth functions $f_r$ for which 
\be\label{actor} \sigma_s (f_r) = e^{2\pi i r\cdot s} \,f_r , \ee and each $f \in 
\cinf(M)$ is the sum of a unique series $f = \sum_{r\in\IZ^n} f_r$, which is 
rapidly convergent in the Fr\'echet topology of $\Cinf(M)$ (see \cite{Rie93a} for 
more details). Let now $\theta = (\theta_{j k} = - \theta_{k j})$ be a real 
antisymmetric $n\times n$ matrix. The $\theta$-deformation of $\cinf(M)$ may be 
defined by replacing the ordinary product by a deformed product, given on 
spectral subspaces by \be \label{eq:star-product} f_r \times_\theta g_{r'} := f_r 
~ \sigma_{\frac{1}{2}r \cdot \theta}( g_{r'}) = e^{ \pi i r \cdot \theta \cdot r' 
} f_r g_{r'} , \ee where $r\cdot \theta$ is the element in $\bR^n$ with 
components $(r\cdot \theta)_k =\sum r_j \theta_{j k}$ for $k=1,\ldots,n$. The 
product in \eqref{eq:star-product} is then extended linearly to all functions in 
$\Cinf(M)$. We denote the space $\cinf(M)$ endowed with the product 
$\times_\theta$ by $\cinf(M_\theta)$. The action $\sigma$ of $\IT^n$ on 
$\cinf(M)$ extends to an action on $\cinf(M_\theta)$ given again by \eqref{actor} 
on the homogeneous elements.

\bigskip Next, let us take $M$ to be a spin manifold with $\cH:=L^2(M,\cS)$ the 
Hilbert space of spinors and $D$ the usual Dirac operator of the metric of $M$. 
Smooth functions act on spinors by pointwise multiplication thus giving a 
representation $\pi : \Cinf(M) \to \B(\cH)$, the latter being the algebra of 
bounded operators on $\ch$.

There is a double cover $c: \widetilde{\bT}^n \to \bT^n$ and a representation of 
$\tilde \bT^n$ on $\cH$ by unitary operators $U(s), s \in \tilde\bT^n$, so that 
\be U(s) D U(s)^{-1} = D, \ee since the torus action is assumed to be isometric, 
and such that for all $f \in \Cinf(M)$, \be U(s) \pi(f) U(s)^{-1} = 
\pi(\sigma_{c(s)}(f)). \ee Recall that an element $T\in \B(\cH)$ is called smooth 
for the action of $\tilde \bT^n$ if the map \bd \tilde \bT^n \ni s \mapsto 
\alpha_s(T) := U(s) T U(s)^{-1}, \ed is smooth for the norm topology. From its 
very definition, $\alpha_s$ coincides on $\pi(C^\infty(M)) \subset \B(\cH)$ with 
the automorphism $\sigma_{c(s)}$. Much as it was done before for the smooth 
functions, we shall use the torus action to give a spectral decomposition of 
smooth elements of $\B(\cH)$.  Any such a smooth element $T$ is written as a 
(rapidly convergent) series $T =\sum T_{r}$ with $r\in\IZ^n$ and each $T_{r}$ is 
homogeneous of degree $r$ under the action of $\tilde \bT^n$, {i.e.} 
\be\label{homocomp} \alpha_s(T_{r}) =e^{2 \pi i r \cdot s } T_{r} , \quad \forall 
~ s \in \tilde \bT^n . \ee Let ($P_1, P_2,\ldots, P_n$) be the infinitesimal 
generators of the action of $\tilde \bT^n$ so that we can write $U(s)=\exp{2 \pi 
i s \cdot P}$.  Now, with $\theta$ a real $n\times n$ anti-symmetric matrix as 
above, one defines a twisted representation of the smooth elements of $\B(\cH)$ 
on $\cH$ by \be\label{twist} L_\theta(T):=\sum_r T_r U( \tfrac{1}{2} r\cdot 
\theta ) = \sum_r T_r \exp \big\{\pi i \, r_j \theta_{jk} P_k \big\}, \ee 
The twist $L_\theta$ commutes with the action $\alpha_s$ of $\tilde \bT^n$ 
and preserves the spectral components of smooth operators: 
\be\label{spectwist} 
\alpha_s(L_\theta(T_r))  = U(s) ~T U(\half r\cdot \theta)~ U(s)^{-1}
=U(s) T U(s)^{-1} U(\half r\cdot \theta) = e^{2 \pi i r \cdot s } L_\theta (T_r).
\ee
Taking smooth functions on $M$ as elements of $\B(\cH)$, via the representation $\pi$, 
the previous definition gives an algebra $L_\theta(\Cinf(M))$ which we may think 
of as a representation (as bounded operators on $\ch$) of the algebra 
$\Cinf(\M)$. Indeed, by the very definition of the product $\times_\theta$ in 
\eqref{eq:star-product} one establishes that \be L_\theta(f\times_\theta g)= 
L_\theta(f) L_\theta(g), \ee proving that the algebra $\Cinf(M)$ equipped with 
the product $\times_\theta$ is isomorphic to the algebra $L_\theta(\Cinf(M))$. It 
is shown in \cite{Rie93a} that there is a natural completion of the algebra 
$\cinf(M_\theta)$ to a $C^*$-algebra $C(M_\theta)$ whose smooth subalgebra -- 
under the extended action of $\IT^n$ -- is precisely $\cinf(M_\theta)$. Thus, we 
can understand $L_\theta$ as a {\it quantization map} from \be L_\theta : 
\Cinf(M) \to \Cinf(\M) , \ee which provides a strict deformation quantization in 
the sense of Rieffel. More generally, in \cite{Rie93a} one considers a (not 
necessarily commutative) $C^*$-algebra $A$ carrying an action of $\bR^n$. For an 
anti-symmetric $n \times n$ matrix $\theta$, one defines a star product 
$\times_\theta$ between elements in $A$ much as we did before. The algebra $A$ 
equipped with the product $\times_\theta$ gives rise to a $C^*$-algebra denoted 
by $A_\theta$. Then the collection $\{A_{\hbar\theta} \}_{\hbar\in [0,1]}$ is a 
continuous family of $C^*$-algebras providing a strict deformation quantization 
in the direction of the Poisson structure on $A$ defined by the matrix $\theta$.

Our case of interest corresponds to the choice $A=C(M)$ with an action of $\bR^n$ 
that is periodic or, in other words, an action of $\bT^n$. The smooth elements in 
the deformed algebra make up the algebra $\Cinf(\M)$. The quantization map will 
play a key role in what follows, allowing us to extend differential geometric 
techniques from $M$ to the noncommutative space $\M$.

\bigskip It was shown in \cite{CL01} that the datum $(L_\theta(C^\infty(M)), \cH, 
D)$ satisfies all properties of a noncommutative spin geometry \cite{C96} (see 
also \cite{GVF01}); there is also a grading $\gamma$ (for the even case) and a 
real structure $J$. In particular, boundedness of the commutators 
$[D,L_\theta(f)]$ for $f \in C^\infty(M)$ follows from 
$[D,L_\theta(f)]=L_\theta([D,f])$, $D$ being of degree $0$ (since $\bT^n$ acts by 
isometries, each $P_k$ commutes with $D$). This noncommutative geometry is an 
isospectral deformation of the classical Riemannian geometry of $M$, in that the 
spectrum of the operator $D$ coincides with that of the Dirac operator $D$ on 
$M$. Thus all spectral properties are unchanged. In particular, the triple is 
$m^+$-summable and there is a noncommutative integral as a Dixmier trace 
\cite{Dix66}, \be\label{dix} \ncint L_\theta (f) := \dix \big( L_\theta(f) 
|D|^{-m} \big), \ee with $f \in \Cinf(\M)$ understood in its representation on 
$\cH$. A drastic simplification of this noncommutative integral is given by the 
Lemma~\cite[Prop. 5.1]{GIV05}. \begin{lma} \label{lma:dix} If $f \in \Cinf(\M)$ 
then $$\ncint L_\theta (f) = \int_{M} f \dd \nu.$$ \end{lma} \begin{proof} Any 
element $f \in \Cinf(\M)$ is given as an infinite sum of functions that are 
homogeneous under the action of $\bT^n$. Let us therefore assume that $f$ is 
homogeneous of degree $k$ so that $\sigma_s(L_\theta(f))=L_\theta(\sigma_s(f))= 
e^{2\pi i k \cdot s} L_\theta (f)$. From the tracial property of the 
noncommutative integral and the invariance of $D$ under the action of $\bT^n$, we 
see that $$ \dix\big(\sigma_s(L_\theta(f))|D|^{-m}\big) = \dix\big(U(s) 
L_\theta(f) U(s)^{-1} |D|^{-m}\big) = \dix (L_\theta(f) |D|^{-m}). $$ In other 
words, $e^{2\pi i k\cdot s} \dix(L_\theta(f) |D|^{-m} )=\dix(L_\theta(f) 
|D|^{-m})$ from which we infer that this trace vanishes if $k\neq 0$. If $k=0$, 
then $L_\theta (f) = f$, leading to the desired result. \end{proof}

\subsection{The manifold $\M$ as a fixed point algebra} A different but 
equivalent approach to these noncommutative manifolds $\M$ was introduced in 
\cite{CD02}. In there the algebra $\Cinf(\M)$ is identified as a fixed point 
subalgebra of $\Cinf(M) \otimes \Cinf(\bT^n_\theta)$ where $\Cinf(\bT^n_\theta)$ 
is the algebra of smooth functions on the noncommutative torus. This 
identification was shown to be useful in extending techniques from commutative 
differential geometry on $M$ to the noncommutative space $\M$.

We recall the definition of the noncommutative $n$-torus $\T^n$ \cite{Rie90a}. 
Let $\theta = (\theta_{jk})$ be a real $n\times n$ anti-symmetric matrix as 
before, and let $\lambda^{jk}=e^{2\pi i\theta_{jk}}$. The unital $*$-algebra 
$\A(\T^n)$ of polynomial functions on $\T^n$ is generated by $n$ unitary elements 
$U^k$, $k=1,\ldots,n$, with relations \be U^j U^k= \lambda^{jk} U^k U^j, \quad 
j,k=1,\ldots,n . \ee The polynomial algebra is extended to the universal 
$C^\ast$-algebra with the same generators. There is a natural action of $\bT^n$ 
on $\A(\T^n)$ by $*$-automorphisms given by $\tau_s (U^k)= e^{2 \pi i s_k } U^k$ 
with $s=(s_k) \in \bT^n$. The corresponding infinitesimal generators $X_k$ of the 
action are algebra derivations given explicitly on the generators by $X_k(U^j)= 
2\pi i \delta_k^j$. They are used \cite{C80} to construct the 
pre-$C^\ast$-algebra $\Cinf(\T^n)$ of smooth functions on $\T^n$, which is the 
completion of $\A(\T^n)$ with respect to the locally convex topology generated by 
the seminorms, \be |u|_r:= \sup_{r_1+\cdots + r_n \leq r}\| X_1^{r_1} \cdots 
X_n^{r_n} (u)\| , \ee and $\| \cdot \|$ is the $C^\ast$-norm. The algebra 
$\Cinf(\T^n)$ turns out to be a nuclear Fr\'echet space and one can unambiguously 
take the completed tensor product $\Cinf(M) \bar\ot \Cinf(\T^n)$. Then, one 
defines $\big(\Cinf(M) \bar\ot \Cinf(\T^n) \big)^{\sigma \ot \tau^{-1}}$ as the 
fixed point subalgebra of $\Cinf(M) \bar\ot \Cinf(\T^n)$ consisting of elements 
$a$ in the tensor product that are invariant under the diagonal action of 
$\bT^n$, {i.e.} such that $\sigma_s \ot \tau_{-s}(a) = a$ for all $s \in \bT^n$. 
The noncommutative manifold $\M$ is defined by ``duality'' by setting \bd 
\Cinf(\M):= \big(\Cinf(M) \bar\ot \Cinf(\T^n) \big)^{\sigma \ot \tau^{-1}} . \ed 
As the notation suggests, the algebra $\Cinf(\M)$ turns out to be isomorphic to 
the algebra $L_\theta (\Cinf(M))$ defined in the previous section.

Next, let $\cS$ be a spin bundle over $M$ and $D$ the Dirac operator on 
$\Gamma^\infty(M,\cS)$, the $\Cinf(M)$-module of smooth sections of $\cS$. The 
action of the group $\bT^n$ on $M$ does not lift directly to the spinor bundle. 
Rather, there is a double cover $c: \widetilde{\bT}^n \to \bT^n$ and a group 
homomorphism $\tilde{s} \to V_{\tilde{s}}$ of $\widetilde{\bT}^n$ into 
$\Aut(\cS)$ covering the action of $\bT^n$ on $M$, \be \label{eq:action-cover} 
V_{\tilde{s}} (f \psi) = \sigma_{c(s)}(f) V_{\tilde{s}}(\psi), \ee for $f \in 
C^\infty(M)$ and $\psi \in \Gamma^\infty(M,\cS)$. According to \cite{CD02}, the 
proper notion of smooth sections $\Gamma^\infty(\M,\cS)$ of a spinor bundle on 
$\M$ are elements of $\Gamma^\infty(M,\cS) \widehat\otimes 
C^\infty(\bT_{\theta/2}^n)$ which are invariant under the diagonal action $V 
\times \tilde{\tau}^{-1}$ of $\tilde{\bT}^n$. Here $\tilde{s} \mapsto 
\tilde\tau_{\tilde{s}}$ is the canonical action of $\tilde{\bT}^n$ on 
$\A(\bT_{\theta/2}^n)$. Since the Dirac operator $D$ commutes with 
$V_{\tilde{s}}$ (remember that the torus action is isometric) one can restrict 
$D\otimes \id$ to the fixed point elements $\Gamma^\infty(\M,\cS)$.

Then, let $L^2(M,\cS)$ be the space of square integrable spinors on $M$ and let 
$L^2(\bT^n_{\theta/2})$ be the completion of $C^\infty(\bT^n_{\theta/2})$ in the 
norm $a \mapsto \|a\|=\tr(a^*a)^{1/2}$, with $\tr$ the usual trace on 
$C^\infty(\bT^n_{\theta/2})$. The diagonal action $V \times \tilde{\tau}^{-1}$ of 
$\tilde{\bT}^n$ extends to $L^2(M,\cS)\otimes L^2(\bT^n_{\theta/2})$ (where it 
becomes $U \times \tau$)  and one defines $L^2(\M, \cS)$ to be the fixed point 
Hilbert subspace. If $D$ also denotes the closure of the Dirac operator on 
$L^2(M,\cS)$, one still denotes by $D$ the operator $D \otimes \id$ on 
$L^2(M,\cS)\otimes L^2(\bT^n_{\theta/2})$ when restricted to $L^2(\M, \cS)$. The 
triple $(\Cinf(\M), L^2(\M,\cS),D)$ is an $m^+$-summable noncommutative spin 
geometry.

\subsection{Vector bundles on $\M$} \label{section:nc-vb} Noncommutative vector 
bundles on $\M$, {i.e.} finite projective modules over $\Cinf(\M)$, were obtained 
in \cite{CD02} as fixed point submodules of $\Gamma^\infty(M,E) \otimes 
\Cinf(\T^n)$ under some diagonal action of the torus $\bT^n$. Here 
$\Gamma^\infty(M,E)$ denotes the $C^\infty(M)$-bimodule of smooth sections of a 
vector bundle $E \to M$. We will presently give an equivalent description of 
these modules over $\Cinf(\M)$ in terms of a kind of $\ast$-product.

Let $E$ be a {\it $\sigma$-equivariant} vector bundle over $M$, that is a bundle 
which carries an action $V$ of $\bT^n$ by automorphisms, covering the action {\it 
$\sigma$} of $\bT^n$ on $M$, \be \label{eq:sigma-equivariant} V_s (f \psi) = 
\sigma_s(f) V_s(\psi), \quad \forall \, f \in \Cinf(M) , \, \psi \in 
\Gamma^\infty(M,E) . \ee We also assume that the topology on the Fr\'echet 
$\Cinf(M)$-bimodule $\Gamma^\infty(M,E)$ is given in terms of a family of 
$\bT^n$-invariant seminorms $\|\cdot \|_i$.

We define the $\Cinf(\M)$-bimodule $\Gamma^\infty(\M, E)$ as the vector space 
$\Gamma^\infty(M,E)$ but with deformed bimodule structure given by 
\begin{align}\label{eq:left-action} f \lt_\theta \psi &= \sum_k f_k 
V_{\frac{1}{2}k\cdot \theta} (\psi) , \\ \psi \rt_\theta f &= \sum_k 
V_{-\frac{1}{2}k\cdot \theta} (\psi) f_k ; \label{eq:right-action} \end{align} 
here $f=\sum_k f_k$ with $f_k \in \Cinf(M)$ homogeneous of degree $k$ under the 
action of $\bT^n$ -- as in \eqref{homocomp} -- and $\psi$ is a smooth section of 
$E$. By using the explicit expression \eqref{eq:star-product} for the star 
product and equation~\eqref{eq:sigma-equivariant}, one checks that these are 
indeed actions of $\Cinf(\M)$. Moreover, the invariance of $\| \cdot \|_i$ under 
the action of $\bT^n$ makes both actions continuous, turning 
$\Gamma^\infty(\M,E)$ into a Fr\'echet $\Cinf(\M)$-bimodule.

The $\Cinf(\M)$-bimodule $\Gamma^\infty(\M,E)$ is finite projective \cite{CD02} 
and still carries an action $V$ of $\bT^n$ with equivariance as in 
\eqref{eq:sigma-equivariant} for both the left and right action of $\Cinf(\M)$: 
the group $\bT^n$ being abelian, one directly establishes that \be V_s (f 
\lt_\theta \psi) = \sigma_s(f) \lt_\theta V_s(\psi), \quad \forall \, f \in 
\Cinf(\M) , \, \psi \in \Gamma^\infty(\M,E), \ee and a similar property for the 
right structure $\rt_\theta$. In fact, since the category of $\sigma$-equivariant 
finite projective module over $\Cinf(\M)$ is equivalent to the category of finite 
projectve modules over $\Cinf(\M) \rtimes_\sigma \bT^n$ (see \cite{Jul81}) for 
all $\theta$ (in particular $\theta=0$), the isomorphism $\Cinf(\M) 
\rtimes_\sigma \bT^n \isom \Cinf(M) \rtimes_\sigma \bT^n$ shows that all 
$\sigma$-equivariant finite projective modules over $\Cinf(\M)$ are of the above 
type \cite[Proof of Prop. 5]{CD02}. This also reflects the result in 
\cite{Rie93b} that the K-groups of a $C^*$-algebra deformed by an action of 
$\bR^n$ are isomorphic to the K-groups of the original $C^*$-algebra: as 
mentioned above, the noncommutative manifolds $\M$ are a special case -- in which 
the starting algebra is commutative and the action periodic -- of the general 
formulation in \cite{Rie93a} of deformations of $C^*$-algebras under an action of 
$\bR^n$.

Although we defined the above left and right actions on sections with respect to 
an action of $\bT^n$ on the vector bundle $E$, the same construction can be done 
for vector bundles carrying an action of the double cover $\tilde \bT^n$. We have 
already seen an example of this double cover action for the spinor bundle, for 
which we defined a left action of $\Cinf(\M)$ using the twisted representation 
\eqref{twist}.

\bigskip From the very definition of $\Gamma^\infty(\M,E)$ the following lemma is 
true. \begin{lma} \label{lma:modulesA} If $E \isom F$ as $\sigma$-equivariant 
vector bundles, then $\Gamma^\infty(\M,E) \isom \Gamma^\infty(\M,F)$ as 
$\Cinf(\M)$-bimodules. \end{lma}

\subsection{Differential calculus on $\M$} \label{section:diff-calc}

It is straightforward to construct a differential calculus on $\M$. This can be 
done in two equivalent manners, either by extending to forms the quantization 
maps, or by using the general construction in \cite{C94} by means of the Dirac 
operator.

Firstly, let $(\Omega(M),\dd)$ be the usual differential calculus on $M$, with 
$\dd$ the exterior derivative. This becomes a Fr\'echet algebra if we consider 
$\Omega(M)$ as the space of smooth sections of a bundle over $M$. Moreover, it 
carries an action $\sigma$ of $\bT^n$ by automorphisms which commutes with the 
differential $\dd$. Again, we assume that the seminorms defining the Fr\'echet 
topology of $\Omega(M)$ are $\bT^n$-invariant.

The quantization map $L_\theta: \Cinf(M) \to \Cinf(\M)$ is extended to 
$\Omega(M)$ by imposing that it commutes with $\dd$. The image 
$L_\theta(\Omega(M))$ will be denoted $\Omega(\M)$ and becomes a Fr\'echet 
algebra with the induced seminorms from $\Omega(M)$. Equivalently, $\Omega(\M)$ 
could be defined to be $\Omega(M)$ as a vector space but equipped with an 
``exterior star product" which is the extension of the product 
\eqref{eq:star-product} to $\Omega(M)$ by the requirement that it commutes with 
$\dd$. Indeed, since the action of $\bT^n$ commutes with $\dd$, an element in 
$\Omega(M)$ can be decomposed into a sum of a rapidly convergent series of 
homogeneous elements for the action of $\bT^n$ -- as was done for $\Cinf(M)$. 
Then one defines a star product $\times_\theta$ on homogeneous elements in 
$\Omega(M)$ as in \eqref{eq:star-product} and denotes $\Omega(\M)=(\Omega(M), 
\times_\theta)$. This construction is in concordance with the previous section, 
when $\Omega(M)$ is considered as a $\Cinf(M)$-bimodule of sections. The extended 
action of $\bT^n$ from $\Cinf(M)$ to $\Omega(M)$ is used to endow the space 
$\Omega(\M)$ with the structure of a $\Cinf(\M)$-bimodule with a left and right 
action given in \eqref{eq:left-action}-\eqref{eq:right-action}.

As mentioned, a differential calculus $\Omega_D(\Cinf(\M))$ on $\Cinf(\M)$ can 
also be obtained from the general procedure \cite{C94} using the isospectral 
Dirac operator $D$ on $\M$ defined above. The $\Cinf(\M)$-bimodule 
$\Omega^p_D(\Cinf(\M))$ of $p$-forms is made of classes of operators \be \omega = 
\sum_j a_0^j [D, a_1^j]\cdots [D, a_p^j], \quad a_i^j \in \Cinf(\M), \ee modulo 
the sub-bimodule of operators \be \big\{ \sum_j [D,b_0^j] [D, b_1^j]\cdots [D, 
b_{p-1}^j]:~ b_i^j \in \Cinf(\M),~ b_0^j [D, b_1^j]\cdots [D, b_{p-1}^j]=0 \big\} 
. \ee The exterior differential $\dd_D$ is given by \be \dd_D\bigg[\sum_j a_0^j 
[D, a_1^j]\cdots [D, a_p^j]\bigg] = \bigg[\sum_j [D,a_0^j] [D, a_1^j]\cdots [D, 
a_p^j]\bigg], \ee and satisfies $\dd_D^2=0$. One also introduces an inner product 
on forms by declaring that forms of different degree are orthogonal, while for 
two $p$-forms $\omega_1,\omega_2$, the product is \be 
\label{def:innerprod-connes} ( \omega_1,\omega_2 )_D=\ncint \omega_1^* \omega_2. 
\ee Here the noncommutative integral is the natural extension of the one in 
\eqref{dix}, \be\label{dix2} \ncint T := \dix \big( T |D|^{-m} \big), \ee with 
$T$ an element in a suitable class of operators. Not surprisingly, these two 
constructions of forms agree \cite{CD02}, that is, the differential calculi 
$\Omega(\M)$ and $\Omega_D(\Cinf(\M))$ are isomorphic. This allows us in 
particular to integrate forms of top dimension, by defining \be \int_{\M} \omega 
:= \ncint \omega_D, \quad \omega \in \Omega^m(\M), \ee where $\omega_D$ denotes the 
element in $\Omega_D^m(\Cinf(\M))$ corresponding to $\omega$ (replacing every $\dd$ 
in $\omega$ by $\dd_D$). We have the following noncommutative Stokes theorem. 
\begin{lma} \label{lma:stokes} If $\omega \in \Omega^{m-1} (\M)$ then $$ 
\int_{\M} \dd \omega = 0 . $$ \end{lma} \begin{proof} From the definition of the 
noncommutative integral, $$ \int_{\M} \dd \omega = \ncint \dd_D\omega_D = \ncint 
\dd_D L_\theta(\omega_D^\class) , $$ with $\omega_D^\class$ the classical 
counterpart of $\omega$, {i.e.} $\omega=L_\theta (\omega_D^\class)$. At this 
point one remembers that $D$ commutes with $L_\theta$ (see 
Section~\ref{subsection:def-torus}), and realizes that there is an analogue of 
Lemma~\ref{lma:dix} for forms, i.e. $\ncint L_\theta (T) = \int_{M} T$. One 
concludes that the above integral vanishes since it vanishes in the classical 
case. \end{proof}

\bigskip The next ingredient is a Hodge star operator on $\Omega(\M)$. 
Classically, the Hodge star operator is a map $\ast: \Omega^p(M) \to 
\Omega^{m-p}(M)$ depending only on the conformal class of the metric on $M$. On 
the one end, since $\bT^n$ acts by isometries, it leaves the conformal structure 
invariant and therefore, it commutes with $\ast$. On the other hand, with the 
isospectral deformation one does not change the metric.  Thus it makes sense to 
define a map $\ast_\theta : \Omega^p(\M)\to \Omega^{m-p}(\M)$ by \be \ast_\theta 
L_\theta(\omega) = L_\theta(\ast \omega), \quad \mathrm{for} \quad \omega \in 
\Omega(M) . \ee

With this Hodge operator, there is an alternative definition of an inner product 
on $\Omega(\M)$. Given that $\ast_\theta$ maps $\Omega^p(\M)$ to 
$\Omega^{m-p}(\M)$, we can define for $\alpha,\beta \in \Omega^p(\M)$ \be 
\label{eq:inner-product-forms} ( \alpha,\beta )_2 = \ncint \ast_\theta(\alpha^* 
\ast_\theta \beta), \ee since $\ast_\theta(\alpha^* \ast_\theta \beta)$ is an 
element in $\Cinf(\M)$. \begin{lma} \label{lma:innerprod-forms} Under the 
isomorphism $\Omega_D(\Cinf(\M)) \simeq \Omega(\M)$, the inner product 
$(\cdot,\cdot)_2$ coincides with $(\cdot,\cdot)_D$. \end{lma} \begin{proof} Let 
$\omega_1, \omega_2$ be two forms in $\Omega_D(\Cinf(M))$, so that 
$L_\theta(\omega_i)$ are two generic forms in $\Omega_D(\Cinf(\M)) \isom L_\theta 
(\Omega(M))=\Omega(\M)$. Then, using Lemma~\ref{lma:dix} it follows that \be 
\ncint L_\theta(\omega_1)^* L_\theta(\omega_2) = \ncint L_\theta(\omega_1^* 
\times_\theta \omega_2) = \ncint \omega_1^* \times_\theta \omega_2, \ee
 Now, the inner product $(~,~)_D$ coincides with $(~,~)_2$ as defined by 
\eqref{eq:inner-product-forms} in the classical case -- under the above 
isomorphism $\Omega_D(\Cinf(M)) \isom \Omega(M)$; see for example 
\cite[VI.1]{C94}.  It follows that the above expression equals \be \ncint 
\ast(\omega_1^* \times_\theta (\ast \omega_2) ) = \ncint 
\ast_\theta(L_\theta(\omega_1)^* (\ast_\theta L_\theta(\omega_2) ) ), \ee using 
Lemma~\ref{lma:dix} for forms once more, together with the defining property of 
$\ast_\theta$. \end{proof}

\begin{lma} \label{lma:d-dstar} The formal adjoint $\dd^*$ of $\dd$ with respect 
to the inner product $(\cdot,\cdot)_2$ -- {i.e.} so that $(\dd^* \alpha, 
\beta)_2=(\alpha,\dd \beta)_2$ -- is given on $\Omega^p(\M)$ by $$ \dd^* = 
(-1)^{m(p+1)+1} \ast_\theta \dd \ast_\theta . $$ \end{lma} \begin{proof} Just as 
in the classical case, this follows from Stokes Lemma~\ref{lma:stokes}, together 
with the observation that $$ \int_{\M} \omega = \ncint *_\theta \omega, \qquad 
\omega \in \Omega^m(\M), $$ again established from the classical case by means of 
the mentioned analogue of Lemma~\ref{lma:dix} for forms. \end{proof}

\begin{rem} The algebra $\Omega(\M)$ can also be defined as a fixed point algebra 
\cite{CD02}. The action $\sigma$ of $\bT^n$ on $\Omega(M)$ allows one to define 
$\Omega(\M)$ by $\big(\Omega(M)\bar\ot \Cinf(\T^n)\big)^{\sigma \ot \tau^{-1}}$. 
Furthermore, since the exterior derivative $\dd$ on $\Omega(M)$ commutes with the 
action of $\bT^n$, the differential $\dd_\theta$, for the fixed point algebra is 
defined as $\dd_\theta=\dd \otimes \id$. Similarly, the Hodge star operator takes 
the form $\ast_\theta = \ast\otimes \id$ with $\ast$ the classical Hodge 
operator. \end{rem}

\section{Gauge theory on the sphere $\S^4$}\label{section:gts}

We now apply the general scheme of noncommutative gauge field theories -- as 
developed in Section~\ref{section:connections} -- to the case of the $\SU(2)$ 
noncommutative principal bundle $\Sk^7 \to \S^4$ constructed in \cite{LS04}. This 
will also make more explicit all the constructions above. It is worth stressing 
that what follows is valid for more general $\theta$-deformed $G$-principal 
bundle. We will come back to this point later in the paper.

\subsection{The principal fibration $\Sk^7 \to \S^4$} The $\SU(2)$ noncommutative 
principal fibration $\Sk^7 \to \S^4$ is given by an algebra inclusion $\A(\S^4) 
\into \A(\Sk^7)$. The algebra $\A(\S^4)$ of polynomial functions on the sphere 
$\S^4$ is
 generated by elements $z_0=z_0^*, z_j, z_j^*$, $j=1,2$, subject to relations 
\be\label{s4t} z_\mu z_\nu = \lambda_{\mu\nu} z_\nu z_\mu, \quad z_\mu z_\nu^* = 
\lambda_{\nu\mu} z_\nu^* z_\mu, \quad z_\mu^* z_\nu^* = \lambda_{\mu\nu} z_\nu^* 
z_\mu^*, \quad \mu,\nu = 0,1,2 , \ee together with the spherical relation 
$\sum_\mu z_\mu^* z_\mu=1$. Here $\theta$ is a real parameter and the deformation 
parameters are given by $\lambda_{\mu\mu}=1$ and \be \lambda_{1 2} = 
\bar{\lambda}_{2 1} =: \lambda=e^{2\pi i \theta}, \quad \lambda_{j 0} = 
\lambda_{0 j } = 1, \quad j=1,2. \ee For $\theta=0$ one recovers the $*$-algebra 
of complex polynomial functions on the usual $S^4$.

The differential calculus $\Omega(\S^4)$ is generated as a graded differential 
$*$-algebra by the elements $z_\mu, z_\mu^*$ in degree 0 and elements $d z_\mu, d 
z_\mu^*$ in degree 1 satisfying the relations, \be\label{rel:diff} 
\begin{aligned} &dz_\mu dz_\nu+ \lambda_{\mu\nu} dz_\nu dz_\mu =0, \\ &z_\mu 
dz_\nu = \lambda_{\mu\nu} dz_\nu z_\mu, \end{aligned} \qquad \begin{aligned} 
&dz_\mu dz_\nu^* + \lambda_{\nu\mu} dz_\nu^* dz_\mu =0;\\ &z_\mu dz_\nu^* = 
\lambda_{\nu\mu} dz_\nu^* z_\mu, \end{aligned} \ee with $\lambda_{\mu\nu}$ as 
before. There is a unique differential $\dd$ on $\Omega(\S^4)$ such that 
$\dd:z_\mu \mapsto d z_\mu$ and the involution on $\Omega(\S^4)$ is the graded 
extension of $z_\mu \mapsto z_\mu^*$: $(\dd \omega)^*=\dd \omega^*$ and 
$(\omega_1\omega_2)^* = (-1)^{d_1 d_2}\omega_2^* \omega_1^*$ for $\omega_j$ a 
form of degree $d_j$.

\bigskip With $\lambda'_{a b} = e^{2 \pi i \theta'_{ab}}$ and $(\theta'_{ab})$ a 
real antisymmetric matrix,
 the algebra $\A(\Sk^7)$ of polynomial functions on the sphere $\Sk^7$ is 
generated by elements $\psi_a, \psi_a^*$, $a=1,\dots,4$, subject to relations 
\be\label{s7t} \psi_a \psi_b = \lambda'_{a b} \psi_b \psi_a, \quad \psi_a 
\psi_b^* = \lambda'_{b a} \psi_b^* \psi_a, \quad \psi_a^*\psi_b^* = \lambda'_{a 
b} \psi_b^* \psi_a^* , \ee and with the spherical relation $\sum_a \psi_a^* 
\psi_a=1$. Clearly, $\A(\Sk^7)$ is a deformation of the $*$-algebra of
 complex polynomial functions on the sphere $S^7$. As before, a differential 
calculus $\Omega(\Sk^7)$ can be defined to be generated by the elements $\psi_a, 
\psi_a^*$ in degree 0 and elements $d \psi_a, d \psi_a^*$ in degree 1 satisfying 
relations similar to the ones in \eqref{rel:diff}.

\bigskip In order to construct the noncommutative fibration over the given 
4-sphere $\S^4$ we need to select
 a particular noncommutative 7 dimensional sphere $\Sk^7$. We take the one 
corresponding to the following deformation parameters \be \label{def:lambda} 
\lambda'_{ab}= \begin{pmatrix} 1 & 1 & \bar\mu & \mu \\ 1 & 1 & \mu & \bar\mu \\ 
\mu &\bar\mu &1 & 1\\ \bar\mu & \mu &1 & 1 \end{pmatrix}, \quad \mu = 
\sqrt{\lambda}, \qquad \mathrm{or} \qquad 
\theta'_{ab}=\frac{\theta}{2}\begin{pmatrix} 0 & 0 & -1 & 1 \\ 0 & 0 & 1 & -1 \\ 
1 & -1 & 0 & 0 \\ -1 & 1 & 0 & 0 \end{pmatrix}. \ee The previous choice is 
essentially the only one\footnote{Compatibility requires that $\mu^2=\lambda$; we 
drop the case $\mu = -\sqrt{\lambda}$ since its ``classical'' limit would 
correspond to ``anti-commuting'' coordinates.} that allows the algebra 
$\A(\Sk^7)$ to carry an action
 of the group $\SU(2)$ by automorphisms and such that the invariant subalgebra 
coincides with
 $\A(\S^4)$.  The best way to see this is by means of the matrix-valued function 
on $\A(\Sk^7)$ given by \be\label{Psi} \Psi = \begin{pmatrix} \psi_1 & - \psi^*_2 
\\ \psi_2 & \psi^*_1 \\ \psi_3 & -\psi^*_4 \\ \psi_4& \psi^*_3 \end{pmatrix}. \ee 
Then the commutation relation of the algebra $\A(\Sk^7)$ gives $\Psi^\dagger \Psi 
= \II_2 $ and $p = \Psi \Psi^\dagger$ is a projection, $p^2=p=p^\dagger$, with 
entries in $\A(\S^4)$. Indeed, the right action of $\SU(2)$ on $\A(\Sk^7)$ is 
simply given by \be \label{actionSU2} \alpha_w (\Psi) = \Psi w , \qquad w = 
\begin{pmatrix} w_1 & -\bar{w}_2 \\ w_2 & \bar{w}_1 \end{pmatrix} \in \SU(2), \ee 
from which the invariance under the $SU(2)$-action of the entries of $p$ follows 
at once. Explicitly, \be \label{projection1} p= \half \begin{pmatrix} 1+z_0 & 0 & 
z_1 & - \bar{\mu} z_2^* \\ 0 & 1+z_0 & z_2 & \mu z_1^* \\ z_1^*& z_2^* & 1-z_0 & 
0\\ -\mu z_2 & \bar{\mu} z_1 & 0 & 1-z_0 \end{pmatrix}, \ee with the generators 
of $\A(\S^4)$ given by \begin{align} \label{subalgebra} z_0 &= \psi^*_1 \psi_1 + 
\psi^*_2 \psi_2 - \psi^*_3 \psi_3 - \psi^*_4 \psi_4& \nn \\
 &= 2(\psi^*_1 \psi_1 + \psi^*_2 \psi_2) -1 = 1 - 2(\psi^*_3 \psi_3 + \psi^*_4 
\psi_4), \nn \\ z_1 &= 2 (\mu \psi_3^* \psi_1 + \psi^*_2 \psi_4)=2(\psi_1 
\psi_3^* + \psi^*_2 \psi_4) , \nn\\ z_2 &= 2(- \psi^*_1 \psi_4 +\bar\mu \psi_3^* 
\psi_2)=2(- \psi^*_1 \psi_4 + \psi_2 \psi_3^*). \end{align} One straightforwardly 
computes that $z_1^* z_1 + z_1^* z_1 + z_0 ^2 = 1$ and the commutation rules $z_1 
z_2 = \lambda z_2 z_1$, $z_1 z_2^* = \bar{\lambda} z_2^* z_1$, and that $z_0$ is 
central.

The relations \eqref{subalgebra} can also be expressed in the form, \be\label{zg} 
z_\mu =\sum_{ab}\psi^*_a(\gamma_\mu)_{ab}\psi_b,\qquad z_\mu^* = 
\sum_{ab}\psi^*_a(\gamma^*_\mu)_{ab}\psi_b, \ee with $\gamma_\mu$ twisted $4 
\times 4$ Dirac matrices given by \begin{align} \label{eq:dirac} 
\gamma_0=\left(\begin{smallmatrix} 1 &&&\\ &1&&\\&&-1&\\&&&-1 
\end{smallmatrix}\right), \qquad \gamma_1=2\begin{pmatrix} 0 & 
\begin{smallmatrix} 0 & 0 \\ 0 & 1 \end{smallmatrix} \\ \begin{smallmatrix} \mu & 
0 \\ 0 & 0 \end{smallmatrix}& 0 \end{pmatrix}, \qquad \gamma_2=2\begin{pmatrix} 0 
& \begin{smallmatrix} 0 & -1 \\ 0 & 0 \end{smallmatrix} \\ \begin{smallmatrix} 0 
& \bar\mu \\ 0 & 0 \end{smallmatrix}& 0 \end{pmatrix}. \end{align} Note that as 
usual $\gamma_0$ is the grading $$ \gamma_0 = 
-\frac{1}{4}[\gamma_1,\gamma_1^*][\gamma_2,\gamma_2^*]. $$ These matrices satisfy 
twisted Clifford algebra relations \cite{CD02}, \be\label{eq:clifford} \gamma_j 
\gamma_k +\lambda_{jk} \gamma_k\gamma_j=0, \qquad \gamma_j \gamma_k^* 
+\lambda_{kj} \gamma_k^*\gamma_j=4\delta_{jk}; \qquad (j,k =1,2). \ee

\bigskip There are compatible toric actions on $\S^4$ and $\Sk^7$. The torus 
$\bT^2$ acts on $\A(\S^4)$ as, \be\label{eq:act-S4} \sigma_s(z_0, z_1, z_2) = 
(z_0, e^{2\pi i s_1} z_1, e^{2\pi i s_2} z_2), \quad s\in\bT^2 . \ee This action 
is lifted to a double cover action on $\A(\Sk^7)$. The double cover map $p: 
\tilde\bT^2 \to \bT^2$ is given explicitly by $p:(s_1,s_2) \mapsto 
(s_1+s_2,-s_1+s_2)$. Then $\tilde\bT^2$ acts on the $\psi_a$'s as, \be 
\label{eq:lift-S7} \tilde\sigma: \left( \psi_1, \psi_2, \psi_3, \psi_4 \right) 
\mapsto \left(e^{2\pi i s_1}~\psi_1, ~e^{-2\pi i s_1}~\psi_2, ~e^{-2\pi i 
s_2}~\psi_3, ~e^{2\pi i s_2}~\psi_4 \right) \ee Equation~\eqref{subalgebra} shows 
that $\tilde\sigma$ is indeed a lifting to $\Sk^7$ of the action of $\bT^2$ on 
$\S^4$. Clearly, this compatibility is built in the construction of the Hopf 
fibration $\Sk^7 \to \S^4$ as a deformation of the classical Hopf fibration $S^7 
\to S^4$ with respect to an action of $\bT^2$, a fact that also dictated the form 
of the deformation parameter $\lambda'$ in \eqref{def:lambda}. As we shall see, 
the previous double cover of tori comes from a spin cover $\Spin_\theta(5)$ of 
$\SO_\theta(5)$ deforming the usual action of $\Spin(5)$ on $S^7$ as a double 
cover of the action of $\SO(5)$ on $S^4$.

In which sense the algebra inclusion $\A(\S^4) \into \A(\Sk^7)$ is a 
noncommutative principal bundle was explained in \cite{LS04} to which we refer 
for more details. Presently we shall recall the construction of associated 
bundles.

\subsection{Associated bundles} \label{section:associated-modules} We shall work 
in the context of smooth functions on $\Sk^7$ and $\S^4$ as defined in general in 
Section~\ref{section:toric-ncm}. Let $\rho$ be any finite-dimensional 
representation of $\SU(2)$ on the vector space $W$. The space that generalizes to 
the case $\theta \neq 0$, the space of $\SU(2)$-equivariant maps from $S^7$ to 
$W$, is given by, \begin{align}\label{eq-map} \Cinf(\Sk^7) \boxtimes_\rho W := 
\big\{\eta \in \Cinf(\Sk^7) \otimes W: (\alpha_w\otimes \id) (\eta) =(\id \otimes 
\rho(w)^{-1})(\eta), \, \forall ~w \in \SU(2) \big\}, \end{align} where 
$\alpha_w$ is the $\SU(2)$ action given in \eqref{actionSU2}. This space is 
clearly a $\Cinf(\S^4)$-bimodule. We have proved in \cite{LS04} that 
$\Cinf(\Sk^7) \boxtimes_\rho W$ is a finite projective module. It is to be 
thought of as the module of sections of a ``noncommutative vector bundle'' 
associated to the principal bundle via the representation $\rho$. It is worth 
stressing that the choice of a projection for a finite projective module requires 
the choice of one of the two (left or right) module structures. Similarly, the 
definition of a Hermitian structure requires the choice of the left or right 
module structure. In the following, we will always work with the right structure 
for the associated modules. There is a natural (right) Hermitian structure on 
$\Cinf(\Sk^7) \boxtimes_\rho W$, defined in terms of a fixed inner product of $W$ 
as, \be \label{def:hest} \langle \eta, \eta' \rangle := \sum_i \bar \eta_i 
\eta'_i. \ee where we denoted $\eta=\sum_i\eta_i \otimes e_i$ and $\eta'=\sum_i 
\eta'_i \otimes e_i $, given an orthonormal basis $\{e_i, \, i=1, \cdots, \dim 
W\}$ of $W$. One quickly checks that $\langle \eta, \eta' \rangle$ is an element 
in $\Cinf(\S^4)$, and that $\langle \mb , \mb \rangle$ satisfies all conditions 
of a right Hermitian structure.

The bimodules $\Cinf(\Sk^7) \boxtimes_\rho W$ are of the type described in 
Section~\ref{section:nc-vb}. The associated vector bundle $E=S^7 \times_\rho W$ 
on $S^4$ carries an action $V$ of $\tilde\bT^2$ induced from its action on $S^7$, 
which is obviously $\sigma$-equivariant. By the very definition of $\Cinf(\Sk^7)$ 
and of $\Gamma^\infty(\S^4,E)$ in Section~\ref{section:nc-vb}, it follows that 
$\Cinf(\Sk^7) \boxtimes_\rho W \isom \Gamma^\infty(\S^4,E)$. Indeed, from the 
undeformed isomorphism, $\Gamma^\infty(S^4,E)\isom \Cinf(S^7)\boxtimes_\rho W$, 
the quantization map $L_{\theta'}$ of $\Cinf(S^7)$, acting only on the first leg 
of the tensor product, establishes this isomorphism, \be 
\label{eq:quantization-vb} L_{\theta'} : \Cinf(S^7) \boxtimes_\rho W \to 
\Cinf(\Sk^7)\boxtimes_\rho W . \ee The above is well defined since the action of 
$\tilde\bT^2$ commutes with the action of $SU(2)$. Also, it is such that 
$L_{\theta'}(f \lt_{\theta'} \eta)=L_{\theta'}(f) L_{\theta'} 
(\eta)=L_{\theta}(f) L_{\theta'} (\eta)$ for $f \in \Cinf(S^4)$ and $\eta \in 
\Cinf(S^7) \boxtimes_\rho W$, due to the identity $L_{\theta'}=L_{\theta}$ on 
$\Cinf(S^4)\subset\Cinf(S^7)$. A similar result holds for the action 
$\rt_{\theta'}$.

\begin{prop} \label{rem:associated-modules} The right $\Cinf(\S^4)$-modules 
$\Cinf(\Sk^7) \boxtimes_\rho W$ admits a homogeneous module basis 
\mbox{$\{e_\alpha, \, \alpha = 1, \cdots, N \}$} -- with a suitable $N$ -- that 
is, under the action $V$ of the torus $\tilde \bT^2$, its elements transform as, 
\be V_{s}(e_\alpha) =e^{2\pi i s\cdot r_{\alpha}}{e_\alpha} , \quad s\in \tilde 
\bT^2 . \ee with $r_\alpha \in \IZ^2$ the degree of $e_\alpha$. \end{prop} 
\begin{proof} The vector space $W$ is a direct sum of irreducible representation 
spaces of $\SU(2)$ and the module $\Cinf(\Sk^7) \boxtimes_\rho W$ decomposes 
accordingly. Thus, we can restrict to irreducible representations. The latter are 
labeled by an integer $n$ with $W\simeq\IC^{n+1}$.

Consider first the case $W=\C^2$. Using \cite[Proposition 2]{LS04}, a basis 
$\{e_1, \cdots, e_4\}$ of the right module $\Cinf(\Sk^7) \boxtimes_\rho \C^2$ is 
given by the columns of $\Psi^\dagger$ where $\Psi$ is the matrix in \eqref{Psi}: 
\be e_1 :=\left(\begin{smallmatrix} ~\psi^*_{1} \\ -\psi_{2} \end{smallmatrix} 
\right), \quad e_2 :=\left(\begin{smallmatrix} \psi^*_{2} \\ \psi_{1} 
\end{smallmatrix} \right), \quad e_3 :=\left(\begin{smallmatrix} ~\psi^*_{3} \\ 
-\psi_{4} \end{smallmatrix} \right), \quad e_4 :=\left(\begin{smallmatrix} 
\psi^*_{4} \\ \psi_{3} \end{smallmatrix} \right) . \ee Using the explicit action 
\eqref{eq:lift-S7} it is immediate to compute the corresponding degrees, \be 
r_1=(-1,0), \quad r_2=(1,0), \quad r_3=(0,1), \quad r_4=(0,-1). \ee More 
generally, with $W=\C^{n+1}$ a homogeneous basis $\{e_\alpha, \alpha = 1, \cdots, 
4^n\}$ for the right module $\Cinf(\Sk^7) \boxtimes_\rho \C^{n+1}$ can be 
constructed from the columns of a similar $(n+1)\times 4^n$ matrix 
$\Psi_\n^\dagger$ given in \cite{LS04}.  \end{proof} The above property allows us 
to prove a useful result for the associated modules. \begin{prop} 
\label{prop:tensor-modules} Let $\rho_1$ and $\rho_2$ be two finite dimensional 
representations of $\SU(2)$ on the vector spaces $W_1$ and $W_2$, respectively. 
There is the following isomorphism of right $\Cinf(\S^4)$-modules, $$ 
\left(\Cinf(\Sk^7) \boxtimes_{\rho_1} W_1 \right) \bar\otimes_{\Cinf(\S^4)} 
\left(\Cinf(\Sk^7) \boxtimes_{\rho_2} W_2 \right) \isom \Cinf(\Sk^7) 
\boxtimes_{\rho_1\otimes \rho_2} (W_1 \otimes W_2). $$ \end{prop} \begin{proof} 
Let $\{ e^1_\alpha, \, \alpha=1,\cdots, N_1 \}$ and $\{ e^2_\beta, \, 
\beta=1,\cdots, N_2\}$ be homogeneous bases for the right modules 
$\Cinf(\Sk^7)\boxtimes_{\rho_1} W_1$ and $\Cinf(\Sk^7) \boxtimes_{\rho_2} W_2$ 
respectively. Then, the right module $\Cinf(\Sk^7) \boxtimes_{\rho_1\otimes 
\rho_2} (W_1 \otimes W_2)$ has a homogeneous basis given by $\{ e^1_\alpha 
\otimes e^2_\beta \}$. We define a map $$ \phi: \left(\Cinf(\Sk^7) 
\boxtimes_{\rho_1} W_1 \right) \bar\otimes_{\Cinf(\S^4)} \left(\Cinf(\Sk^7) 
\boxtimes_{\rho_2} W_2 \right) \to \Cinf(\Sk^7) \boxtimes_{\rho_1\otimes \rho_2} 
(W_1 \otimes W_2) $$ by $$ \phi(e^1 _\alpha \times_\theta f^1_\alpha \otimes 
e^2_\beta \times_\theta f^2_\beta) = (e^1_\alpha \otimes e^2_\beta) \times_\theta 
\sigma_{r_\beta \theta}(f^1_\alpha) \times_\theta f^2_\beta, $$ with summation 
over $\alpha$ and $\beta$ understood. Here $r_\beta \in \bZ^2$ is the degree of 
$e^2_\beta$ under the action of $\tilde\bT^2$, so that $e^2_\beta \times_\theta 
\sigma_{r_\beta \theta} (f) = f \times_\theta e^2_\beta$ for any $f \in 
\Cinf(\S^4)$. Note that this map is well-defined since $$ \phi\left(e^1_\alpha 
\times_\theta f^1_\alpha \times_\theta f \otimes e^2_\beta \times_\theta 
f^2_\beta - e^1_\alpha \times_\theta f^1_\alpha \otimes f \times_\theta e^2_\beta 
\times_\theta f^2_\beta\right)=0. $$ Moreover, it is clearly a map of right 
$\Cinf(\S^4)$-modules. In fact, it is an isomorphism with its inverse given 
explicitly by $$ \phi^{-1} \left( (e^1_\alpha \otimes e^2_\beta) \times_\theta 
f_{\alpha\beta} \right) = e_1^1 \otimes (e^2_\beta \times_\theta f_{1\beta}) + 
\cdots + e_{N_1}^1 \otimes (e^2_\beta \times_\theta f_{N_1\beta}) $$ with 
$f_{\alpha\beta} \in \Cinf(\S^4)$. \end{proof}

\noindent Given a right $\Cinf(\S^4)$-module $\E$, its {\it dual module} is 
defined by \be \E' := \big\{ \phi: \E \to \Cinf(\S^4) : \phi(\eta f ) = 
\phi(\eta) f, \, \forall ~f \in\Cinf(\S^4) \big\}, \ee and is naturally a left 
$\Cinf(\S^4)$-module. In the case that $\E$ is also a left $\Cinf(\S^4)$-module, 
then $\E'$ is also a right $\Cinf(\S^4)$-module. If $\E:=\Cinf(\Sk^7) 
\boxtimes_\rho W$ comes from the $\SU(2)$-representation $(W,\rho)$, by using the 
induced dual representation $\rho'$ on the dual vector space $W'$ given by \be 
\big( \rho'(w) v'\big) (v) := v'\big(\rho(w)^{-1} v\big); \qquad \forall ~v' \in 
W', v \in W , \ee we have that \begin{align} \E' &\isom \Cinf(\Sk^7) 
\boxtimes_{\rho'} W' \nn \\ &:= \big\{ \phi \in \Cinf(\Sk^7) \otimes W': 
(\alpha_w \otimes \id) (\phi) =(\id \otimes \rho'(w)^{-1})(\phi), \, \forall ~w 
\in \SU(2) \big\} . \end{align} Next, let $L(W)$ denote the space of linear maps 
on $W$, so that $L(W)=W\otimes W'$. The adjoint action of $\SU(2)$ on $L(W)$ is 
the tensor product representation $\ad:=\rho \otimes \rho'$ on $W \otimes W'$. We 
define \begin{align} \Cinf(\Sk^7) \boxtimes_{\ad} L(W) := \big\{ & T \in 
\Cinf(\Sk^7) \otimes L(W) : \nn \\ & \qquad : (\alpha_w \otimes \id) (T) =(\id 
\otimes \ad(w)^{-1})(T), \, \forall ~w \in \SU(2) \big\} , \end{align} and write 
$T=T_{ij} \otimes e_{ij}$ with respect to the basis $\{ e_{ij} \}$ of $L(W)$ 
induced from the basis $\{e_i\}_{i=1}^{\dim W}$ of $W$ and the dual one 
$\{e'_i\}_{i=1}^{\dim W}$ of $W'$.

On the other hand, there is the endomorphism algebra \be \End_{\Cinf(\S^4)}(\E):= 
\big\{ T: \E \to \E : T(\eta f ) = T (\eta) f, \, \forall ~f \in\Cinf(\S^4) 
\big\}.  \ee We will suppress the subscript $\Cinf(\S^4)$ from $\End$ in the 
following. As a corollary to the previous Proposition, we have the following. 
\begin{prop} \label{prop:end} Let $\E:=\Cinf(\Sk^7) \boxtimes_\rho W$ for a 
finite-dimensional representation $\rho$. Then there is an isomorphism of 
algebras \[ \End(\E) \isom \Cinf(\Sk^7) \boxtimes_{\ad} L(W). \] \end{prop} 
\begin{proof} By Proposition~\ref{prop:tensor-modules}, we have that 
$\Cinf(\Sk^7) \boxtimes_{\ad} L(W) \equiv \Cinf(\Sk^7) \boxtimes_{\rho\otimes 
\rho'} (W \otimes W')$ is isomorphic to $\E \bar\otimes_{\Cinf(\S^4)} \E'$ as a 
right $\Cinf(\S^4)$-module. Since $\E$ is a finite projective 
$\Cinf(\S^4)$-module, there is an isomophism $\End(\E) \isom \E 
\bar\otimes_{\Cinf(\S^4)} \E'$, 
whence the result. 
\end{proof} We see that the algebra of endomorphisms of the module $\E$ can be 
understood as the space of sections of the noncommutative vector bundle 
associated to the adjoint representation on $L(W)$ -- exactly as it happens in 
the classical case. This also allows an identification of skew-hermitian 
endomorphisms $\End^s(\E)$ -- which were defined in general in \eqref{def:skew} 
-- for the toric deformations at hand. \begin{corl} There is an identification
 \[
 \End^s(\E) \isom \Cinf_\bR(\Sk^7) \boxtimes_\ad u(n),
 \]
  with $\Cinf_\bR(\Sk^7)$ denoting the subspace of self-adjoint elements in 
$\Cinf(\Sk^7)$ and $u(n)$ consists of skew-adjoint matrices in $M_n(\C) \isom 
L(W)$, with $n=\dim W$. \end{corl} \begin{proof} Note that the involution $T 
\mapsto T^*$ in $\End(\E)$ reads in components $T_{ij} \mapsto \bar{T_{ji}}$ so 
that, with the identification of Proposition~\ref{prop:end}, the space 
$\End^s(\E)$ is made of elements $X \in \Cinf(\Sk^7) \boxtimes_\ad L(W)$ 
satisfying $\bar{X_{ji}} = -X_{ij}$. Since any element in $\Cinf(\Sk^7)$ can be 
written as the sum of two self-adjoint elements, $X_{ij} = X_{ij}^\Re + i 
X_{ij}^\Im$, we can write $$ X=\sum_i X_{ii}^\Im\otimes i e_{ii}+ \sum_{i \neq j} 
X_{ij}^\Re \otimes (e_{ij} - e_{ji})  + X_{ij}^\Im \otimes (i e_{ij}+i e_{ji}) = 
\sum_a X_a \otimes \sigma^a, $$ with $X_a$ generic elements in $\Cinf_\bR(\Sk^7)$ 
and $\{\sigma^a, a=1,\ldots,n^2\} $ the generators of $u(n)$. \end{proof}

\begin{ex} \label{ex:instanton-bundle}
 Of central interest in the following is the special case of the noncommutative 
instanton bundle first constructed in \cite{CL01}. Now $W=\C^2$ and $\rho$ is the 
defining representation of $\SU(2)$. The projection $p$ giving the 
$\Cinf(\S^4)$-module $\Cinf(\Sk^7)\boxtimes_\rho \C^2$ as a direct summand of 
$\big( \Cinf(\S^4) \big)^N$ for some $N$, is precisely given by the one in 
\eqref{projection1} and $N=4$. Indeed, a generic element in 
$\Cinf(\Sk^7)\boxtimes_\rho \C^2$ is of the form $\Psi^\dagger f$ for some $f \in 
\Cinf(\S^4) \otimes \C^4$ with $\Psi$ defined in \eqref{Psi}, and the 
correspondence is given by \be \Cinf(\Sk^7)\boxtimes_\rho \C^2 \isom p \big( 
\Cinf(\S^4) \big)^4, \quad \Psi^\dagger f \leftrightarrow p f. \ee Furthermore, 
$\End(\E) \isom \Cinf(\Sk^7)\boxtimes_{\ad} M_2(\C)$. It is a known fact that 
$M_2(\C)$ decomposes into the adjoint representation $su(2)$ and the trivial 
representation $\C$ while it is easy to see that $\Cinf(\Sk^7 )\boxtimes_{\id} \C 
\isom \Cinf(\S^4)$. Thus, we conclude that \be \End(\E) \isom 
\Gamma^\infty(\ad(\Sk^7)) \oplus \Cinf(\S^4), \ee where we have set 
$\Gamma^\infty(\ad(\Sk^7)):=\Cinf(\Sk^7)\boxtimes_{\ad} su(2)$. The latter 
$\Cinf(\S^4)$-bimodule will be understood as the space of (complex) sections of 
the adjoint bundle. It is the complexification of the traceless skew-hermitian 
endomorphisms $\Cinf_\bR(\Sk^7) \boxtimes_\ad su(2)$. \end{ex}

\subsection{Yang--Mills theory on $\S^4$} \label{section:ym} Let us now move to 
the main goal of this paper and discuss the Yang--Mills action functional on 
$\S^4$ together with its equations of motion. We will see that instantons 
naturally arise as the local minima of this action.

Before we proceed we recall the noncommutative spin structure 
$(\Cinf(\S^4),\ch,D,\gamma_5$) of $\S^4$ with $\ch=L^2(S^4,\cs)$ the Hilbert 
space of spinors, $D$ the undeformed Dirac operator, and $\gamma_5$ -- the even 
structure -- the fifth Dirac matrix.

Let $\E=\Gamma^\infty(\S^4,E)$ for some $\sigma$-equivariant vector bundle $E$ on 
$S^4$, so that there exists a projection $p \in M_N(\Cinf(\S^4))$ such that 
$\E\isom p (\Cinf(\S^4)^N$. Recall from Section~\ref{section:connections} that a 
connection $\nabla$ on $\E=\Gamma^\infty(\S^4,E)$ can be given as a map from $\E$ 
to $\E\otimes \Omega^1(\S^4)$ satisfying a Leibniz rule with respect to the right 
multiplication of $\Cinf(\S^4)$ on $\ce$. It is also required to be compatible 
with a Hermitian structure on $\ce$. The Yang--Mills action functional is defined 
in terms of the curvature $F$ of $\nabla$, an element in 
$\Hom_{\Cinf(\S^4)} (\E, \E \otimes \Omega^2(\S^4) )$. Equivalently, it can be 
thought of as belonging to $\End_{\Omega(\S^4)}(\E\ot \Omega(\S^4))$ of degree 2. 
We define an inner product on the latter algebra following \cite[III.3]{C94}. Any 
$T \in \End_{\Omega(\S^4)}(\E\ot \Omega(\S^4))$ of degree $k$ can be also 
understood as an element in $p M_N(\Omega^k(\S^4))p$, since $\E\ot \Omega(\S^4)$ 
is a finite projective module over $\Omega(\S^4)$. A trace over internal indices 
together with the inner product in \eqref{eq:inner-product-forms}, defines the 
inner product $(\cdot,\cdot)_2$ on $\End_{\Omega(\S^4)}(\E\ot \Omega(\S^4))$. In 
particular, we can give the following definition: \begin{defn} \label{def:YM} The 
{\rm Yang--Mills action functional} on the collection $C(\ce)$ of compatible 
connections on $\E$ is defined by $$ \YM(\nabla)=\big( F,F\big)_2 = \ncint 
*_\theta \tr (F *_\theta F), $$ for any connection $\nabla$ with curvature $F$. 
\end{defn} \noindent Recall from Section~\ref{section:connections} that gauge 
transformations are unitary endomorphisms $\cU(\E)$ of $\E$.

\begin{lma} The Yang--Mills action functional is gauge invariant, positive and 
quartic. \end{lma} \begin{proof} From equation~\eqref{ugcur}, under a gauge 
transformation $u \in \cU(\E)$ the curvature $F$ transforms as $F \mapsto u^* F 
u$ . Since $\cU(\E)$ can be identified with the unitary elements in $p M_N(\A) 
p$, it follows that \begin{align*} \YM(\nabla^u)=\ncint \sum_{i,j,k,l} *_\theta 
(\bar{u_{ji}} F_{jk} *_\theta F_{kl} u_{li} )=\YM(\nabla), \end{align*} using the 
tracial property of the Dixmier trace and the fact that $u_{li} 
\bar{u_{ji}}=\delta_{lj}$. \\ Positiveness of the Yang--Mills action functional 
follows from Lemma~\ref{lma:innerprod-forms} giving $$ (F,F)_2=(F_D,F_D)_D=\ncint 
F_D^* F_D, $$ which is clearly positive. \end{proof}

\bigskip The Yang--Mills equations of motion (equations for critical points) are 
obtained from the Yang--Mills action functional by a variational principle. Let 
us describe how this principle works in our case. We consider a linear 
perturbation $\nabla_t=\nabla + t \alpha$ of a connection $\nabla$ on $\E$ by an 
element $\alpha \in \Hom(\E,\E\ot_{\Cinf(\S^4)} \Omega^1(\S^4)) $. The curvature 
$F_t$ of $\nabla_t$ is readily computed as $F_t = F + t [\nabla, \alpha] + \cO 
(t^2)$. If we suppose that $\nabla$ is a critical point of the Yang--Mills action 
functional, this linear perturbation should not affect the action. In other 
words, we need \be \frac{\partial}{\partial t}\bigg|_{t=0} \YM(\nabla_t)=0. \ee 
If we substitute the explicit formula for $F_t$, we obtain \be \big( [\nabla, 
\alpha],F \big)_2 + \bar{\big( [\nabla, \alpha] , F\big)_2} = 0, \ee using the 
fact that $(\cdot,\cdot)_2$ defines a complex scalar product on $\Hom(\E,\E\ot 
\Omega(\S^4))$. Positive definiteness of this scalar product implies that 
$(F_t,F_t)=\bar{(F_t,F_t)}$, which when differentiated with respect to $t$, at 
$t=0$, yields $\big( [\nabla, \alpha],F \big)_2 = \bar{\big( [\nabla, \alpha] , 
F\big)_2}$; hence, $( [\nabla, \alpha],F )_2=0$. Using the fact that $\alpha$ was 
arbitrary, we derive the following equations of motion \be [\nabla^*, F \big] = 
0. \ee where the adjoint $[\nabla^*, \cdot ]$ is defined with respect to the 
inner product $(\cdot,\cdot)_2$, {i.e.} \be\label{adjoint-conn} \big( 
[\nabla^*,\alpha], \beta \big)_2 = \big( \alpha, [\nabla, \beta] \big)_2 \ee for 
any $\alpha \in \Hom(\E,\E \ot \Omega^3(\S^4))$ and any $\beta \in \Hom(\E,\E\ot 
\Omega^1(\S^4))$. From Lemma~\ref{lma:d-dstar}, it follows that $[\nabla^* ,F]= 
\ast_\theta [\nabla, \ast_\theta F]$, so that the equations of motion can be 
written as the more familiar \textit{Yang--Mills equations}, \be \label{eq:ym} 
[\nabla, \ast_\theta F]=0. \ee Clearly, connections with a self-dual or 
antiself-dual curvature, $\ast_\theta F= \pm F$, are special solutions of the 
Yang--Mills equations. Indeed, in this case the latter equations follow directly 
from the Bianchi identity, $[\nabla,F]=0$, given in Proposition~\ref{ubianchi}.

\bigskip We will now establish a connection between the Yang--Mills action 
functional and the so-called topological action \cite[VI.3]{C94} on $\S^4$. 
Suppose $\E$ is a finite projective module over $\Cinf(\S^4)$ defined by a 
projection $p \in M_N(\Cinf(\S^4))$. The topological action for $\E$ is the 
pairing between the class of $p$ in K-theory and the cyclic cohomology of 
$\Cinf(\S^4)$. For computational purposes, we define it in terms of the curvature 
of a connection on $\E$, \begin{defn}\label{topact} Let $\nabla$ be a connection 
on $\E$ with curvature $F$. The {\rm topological action} is given by $$ 
\Top(\E)=(F,\ast_\theta F)_2=\ncint \ast_\theta \tr(F^2) $$ where the trace is 
taken over internal indices and in the second equality we have used the identity 
$\ast_\theta \circ \ast_\theta = \id$ on $\S^4$. \end{defn} Let us show that this 
does not depend on the choice of a connection on $\E$. Since two connections 
differ by an element $\alpha$ in $\Hom_{\Cinf(\S^4)}(\E,\E\otimes 
\Omega^1(\S^4))$, we have to establish that $(F',\ast_\theta F')_2=(F,\ast_\theta 
F)_2$ where $F'=F+t [\nabla,\alpha] + \cO(t^2)$ is the curvature of 
$\nabla':=\nabla + t \alpha$, $t \in \bR$ . By definition of the inner product 
$(\cdot,\cdot)_2$ we have that \begin{align*} (F', \ast_\theta F')_2 - (F, 
\ast_\theta F)_2 &= t (F, \ast_\theta [\nabla, \alpha])_2+ 
t([\nabla,\alpha],\ast_\theta F)_2 +\cO(t^2) \nn \\ &= t(F, [\nabla^*,*_\theta 
\alpha])_2+ t([\nabla^*,*_\theta\alpha],F)_2 +\cO(t^2). \end{align*} With 
equation~\eqref{adjoint-conn}, this reduces to \begin{align*} (F', \ast_\theta 
F')_2 - (F, \ast_\theta F)_2 &= t([\nabla,F], *_\theta \alpha)_2+ 
t(*_\theta\alpha,[\nabla,F])_2 +\cO(t^2), \end{align*} which vanishes modulo 
$t^2$ due to the Bianchi identity $[\nabla,F]=0$.

The Hodge star operator $\ast_\theta$ splits $\Omega^2(\S^4)$ into a self-dual 
and antiself-dual space, \be \Omega^2(\S^4)=\Omega_+^2(\S^4) \op 
\Omega_-^2(\S^4). \ee In fact, $\Omega_\pm^2(\S^4)=L_\theta 
\left(\Omega_\pm^2(S^4) \right)$. This decomposition is orthogonal with respect 
to the inner product $(\cdot,\cdot)_2$ a fact that follows from the property 
$(\alpha,\beta)_2=\bar{(\beta,\alpha)_2}$; thus we can write the Yang--Mills 
action functional as \be \YM(\nabla) = \big( F_+,F_+\big)_2 + \big( 
F_-,F_-\big)_2. \ee Comparing this with the analogue decomposition of the 
topological action, \be \Top(\E)=\big( F_+,F_+\big)_2 - \big( F_-,F_-\big)_2, \ee 
we see that $\YM(\nabla) \geq | \Top(\E) |$, with equality holding iff \be 
\ast_\theta F = \pm F. \ee Solutions of these equations are called instantons. We 
conclude that instantons are absolute minima of the Yang--Mills action 
functional.

\section{Construction of $\SU(2)$-instantons on $\S^4$}\label{sect:constrinst}

In this section, we construct a set of $\SU(2)$ instantons on $\S^4$, with 
topological charge $1$, by acting with a twisted infinitesimal conformal symmetry 
on the basic instanton on $\S^4$ constructed in \cite{CL01}. We will find a 
five-parameter family of such instantons. Then we prove that the ``tangent 
space'' to the moduli space of irreducible instantons at the basic instanton is 
five-dimensional, proving that this set is complete. Here, one has to be careful 
with the notion of tangent space to the moduli space. As will be discussed 
elsewhere \cite{LPRS05}, one can construct a noncommutative family of instantons, 
that is instantons parametrized by the quantum quotient space of the deformed conformal group 
$\SL_\theta(2,\bH)$ by the deformed spin group $\Sp_\theta(2)$. It turns out that 
the basic instanton of \cite{CL01} is a ``classical point'' in this moduli space 
of instantons. We perturb the connection $\nabla_0$ linearly by sending $\nabla_0 
\mapsto \nabla_0+ t \alpha$ where $t\in \bR$ and $\alpha \in \Hom(\E,\E \ot 
\Omega^1(\S^4))$. In order for this new connection still to be an instanton, we 
have to impose the self-duality equations on its curvature. After deriving this 
equation with respect to $t$, setting $t=0$ afterwards, we obtain the linearized 
self-duality equations to be fulfilled by $\alpha$. It is in this sense that we 
are considering the tangent space to the moduli space of instantons at the point 
$\nabla_0$.

\subsection{The basic instanton} We start with a technical lemma that simplifies 
the discussion. Let $\E = \Cinf(\Sk^7)\boxtimes_\rho W$ be a module of sections 
associated to a finite dimensional representation of $\SU(2)$, as defined in 
equation~\eqref{eq-map}.
 
\begin{lma} There is the following isomorphism of right $\Cinf(\S^4)$-modules, $$ 
\E \ot_{\Cinf(\S^4)} \Omega(\S^4) \isom \Omega(\S^4)\ot_{\Cinf(\S^4)}\E. $$ 
Consequently, $\Hom(\E,\E\ot_{\Cinf(\S^4)} \Omega(\S^4)) \isom 
\Omega(\S^4)\ot_{\Cinf(\S^4)} \End(\E)$. \end{lma} \begin{proof} Recall from 
Proposition~\ref{rem:associated-modules} that the right $\Cinf(\S^4)$-module $\E$ 
has a homogeneous module-basis $\{e_\alpha,\alpha=1,\cdots, N\}$ for some $N$ and 
each $e_\alpha$ of degree $r_\alpha$. A generic element in $\E \ot_{\Cinf(\S^4)} 
\Omega(\S^4)$ can be written as a sum $\sum_\alpha e_\alpha \otimes_{\Cinf(\S^4)} 
\omega^\alpha$, with $\omega^\alpha$ an element in $\Omega(\S^4)$. Now, for every 
$\omega \in \Omega(\S^4)$ there is an element $\tilde \omega \in \Omega(\S^4)$ -- 
given explicitly by $\tilde\omega = \sigma_{r_\alpha \cdot \theta} (\omega)$ -- 
such that $e_\alpha \omega = \tilde \omega e_\alpha$ where the latter equality 
holds inside the algebra $\Omega(\Sk^7)$ (recall that $\Omega(\S^4) \subset 
\Omega(\Sk^7)$). We can thus define a map $$ T: \E \ot_{\Cinf(\S^4)} \Omega(\S^4) 
\isom \Omega(\S^4)\ot_{\Cinf(\S^4)}\E, $$ by $T(e_\alpha \otimes_{\Cinf(\S^4)} 
\omega^\alpha) = \tilde \omega^\alpha \otimes e_\alpha$; it is a right 
$\Cinf(S^4)$-module map: $$ T\left(e_\alpha \otimes_{\Cinf(\S^4)} (\omega^\alpha 
\times_\theta f) \right) = (\tilde \omega^\alpha \times_\theta \tilde f) 
\otimes_{\Cinf(\S^4)} e_\alpha = T(e_\alpha \otimes_{\Cinf(\S^4)} \omega^\alpha) 
\times_\theta f. $$ Since an inverse map $T^{-1}$ is easily constructed, $T$ 
gives the desired isomorphism. \end{proof} \noindent Thus, we can unambiguously 
use the notation $\Omega(\S^4,\E)$ for the above right $\Cinf(\S^4)$-module $\E 
\ot_{\Cinf(\S^4)} \Omega(\S^4)$.

We let $\nabla_0=p\circ \dd$ be the canonical (Grassmann) connection on the 
projective module $\E=\Cinf(\Sk^7)\boxtimes_\rho \C^2 \simeq p (\Cinf(\S^4) )^4$, 
with the projection $p= \Psi^\dagger \Psi$ of \eqref{projection1} and $\Psi$ is 
the matrix \eqref{Psi} (refer to Example~\ref{ex:instanton-bundle}). When acting 
on equivariant maps, we write $\nabla_0$ as \be\label{cancon} \nabla_0: \E \to \E 
\ot_{\Cinf(\S^4)} \Omega^1(\S^4), \qquad (\nabla_0 f)_i = \dd f_i + \omega_{ij} 
\times_\theta f_j, \ee where $\omega$ -- referred to as the gauge potential -- is 
given in terms of the matrix $\Psi$ by \be\label{cangp} \omega=\Psi^\dagger \dd 
\Psi . \ee The above, is a $2\times 2$-matrix with entries in $\Omega^1(\Sk^7)$ 
satisfying $\bar{\omega_{ij} }=\omega_{ji}$ and $\sum_i \omega_{ii}=0$. Note here 
that the entries $\omega_{ij}$ commute with all elements in $\Cinf(\Sk^7)$. 
Indeed, from \eqref{Psi} we see that the elements in $\omega_{ij}$ are 
$\bT^2$-invariant and hence central (as one forms) in $\Omega(\Sk^7)$. In other 
words $L_\theta(\omega)=\omega$, which shows that for an element $f \in \E$ as 
above, we have $\nabla_0(f)_i=\dd f_i + \omega_{ij} \times_\theta f_j= \dd f_i + 
\omega_{ij} f_j$ which coincides with the action of the classical connection 
$\dd+\omega$ on $f$.

The curvature $F_0=\nabla_0^2=\dd \omega + \omega^2$ of $\nabla_0$ is an element 
of $\End(\E) \ot_{\Cinf(\S^4)} \Omega^2(\S^4)$ that satisfies \cite{AB02,CD02} 
the self-duality equation \be \ast_\theta F_0 = F_0; \ee hence this connection is 
an instanton. At the classical value of the deformation parameter, $\theta=0$, 
the connection \eqref{cangp} is nothing but the $SU(2)$ instanton of \cite{BPST}.

Its ``topological charge'', i.e. the values of $\Top(\E)$ in 
Definition~\ref{topact}, was already computed in \cite{CL01}. Clearly it depends 
only on the class $[p]$ of the bundle and can be evaluated as the index 
\be\label{topbasic} \Top([p]) = \ind(D_p) = \ncint \gamma_5 \pi_D(\chern_2(p)) 
\ee
 of the twisted Dirac operator \[ D_p = p (D \ot \II_4) p . \] The last equality 
in \eqref{topbasic} follows from the vanishing of the class $\chern_1(p)$ of the 
bundle\footnote{The Chern character classes and their realization as operators 
are in the Appendix.}. On the other hand, one finds \[ \pi_D 
\big(\chern_2(p))\big) = 3 \gamma_5, \] which, together with the fact that \[ 
\ncint 1 = \dix |D|^{-4} = \frac{1}{3}, \] on $S^4$ (see for instance 
\cite{GVF01,Lnd97}), gives the value $\Top([p])=1$.

We aim at constructing all connections $\nabla$ on $\E$ whose curvature satisfies 
the self-duality equations and have topological charge equal to $1$. We can write 
any such connection in terms of the canonical connection as in 
equation~\eqref{uconn}, {i.e.} $\nabla=\nabla_0 + \alpha$ with $\alpha$ a 
one-form valued endomorphism of $\E$. Clearly, this will not change the value of 
the topological charge. Being interested in $\SU(2)$-instantons, we impose that 
$\alpha$ is traceless and skew-hermitian. Here the trace is taken in the second 
leg of $\End (\E) \isom \Cinf(\Sk^7) \boxtimes_\ad M_2(\C)$. When complexified, 
we get an element $\alpha \in 
\Omega^1(\S^4)\ot_{\Cinf(\S^4)}\Gamma^\infty(\ad(\Sk^7)) =: \Omega^1(\ad(\S^4))$ 
(cf. Example~\ref{ex:instanton-bundle}).

As usual, we impose an irreducibility condition on the instanton connections, a 
connection on $\E$ being {\it irreducible} if it cannot be written as the sum of 
two other connections on $\E$. We are interested only in the irreducible 
instanton connections on the module $\E$.

\begin{rem}\label{rem:conn-form-central-general} In \cite{LS04}, we constructed 
projections $p_{(n)}$ for all modules $\Cinf(\Sk^7) \boxtimes_{\rho} \C^n$ over 
$\Cinf(\S^4)$ associated to the irreducible representations $\C^n$ of $\SU(2)$. 
The induced Grassmann connections $\nabla_0^{(n)}:= p_{(n)} \dd$, when acting on 
$\Cinf(\Sk^7) \boxtimes_{\rho} \C^n$, were written as $\dd + \omega_{(n)}$, with 
$\omega_{(n)}$ an $n\times n$ matrix with entries in $\Omega^1(\Sk^7)$. A similar 
argument as above then shows that all $\omega_{(n)}$ have entries that are 
central (as one forms) in $\Omega(\Sk^7)$; again, this means that 
$L_\theta(\omega_\n)=\omega_\n$. In particular, this holds for the adjoint bundle 
associated to the adjoint representation on $su(2)_\C \isom \C^3$ (as complex 
representation spaces), from which we conclude that $\nabla_0^{(2)}$ coincides 
with $[\nabla_0,\cdot]$ (since this is the case if $\theta=0$). \end{rem}

\subsection{Twisted infinitesimal symmetries} The noncommutative sphere $\S^4$ 
can be realized as a quantum homogeneous space of the quantum orthogonal group 
$\SO_\theta(5)$ \cite{Var01,CD02,AB02}. In other words, $\A(\S^4)$ can be 
obtained as the subalgebra of $\A(\SO_\theta(5))$ made of elements that are 
coinvariant under the natural coaction of $\SO_\theta(4)$ on $\SO_\theta(5)$. For 
our purposes, it turns out to be more convenient to take a dual point of view and 
consider an {\it action} instead of a coaction. We obtain a twisted symmetry 
action of the Lie algebra $so(5)$ on $\S^4$. Elements of $so(5)$ act as twisted 
derivations on the algebra $\A(\S^4)$.  This action is lifted to $\Sk^7$ and the 
basic instanton $\nabla_0$ described above is invariant under this infinitesimal 
twisted symmetry.

Different instantons are obtained by a twisted symmetry action of $so(5,1)$. 
Classically, $so(5,1)$ is the conformal Lie algebra consisting of the 
infinitesimal diffeomorphisms leaving the conformal structure invariant. The Lie 
algebra $so(5,1)$ is given by adding 5 generators to $so(5)$. We explicitly 
describe its action on $\S^4$ together with its lift to $\Sk^7$ as an algebra of 
twisted derivations. The induced action on forms leaves the conformal structure 
invariant and when acting on $\nabla_0$ eventually results in a five-parameter 
family of instantons.

In fact, what we are really describing are Hopf algebras $\U_\theta(so(5))$ and 
$\U_\theta(so(5,1))$ which are obtained from the undeformed Hopf algebras 
$\U(so(5))$ and $\U(so(5,1))$ via a twist of a Drinfel'd type. Twisting of 
algebras and coalgebras has been known for some time \cite{Dr83,Dr90,GZ94}. The 
twists relevant for toric noncommutative manifolds are associated to the Cartan 
subalgebra of a Lie algebra and were already introduced in \cite{Re90}. Their use 
to implement symmetries of toric noncommutative manifolds like the ones of the 
present paper was made explicit in \cite{Sit01}. The geometry of multi-parametric 
quantum groups and quantum enveloping algebras coming from twists has been 
studied in \cite{AC1,AC2}.

Since we do not explicitly need Hopf algebras, we postpone a full fledged 
analysis of them and of their actions which will appear elsewhere. Here we shall 
rather give the explicit actions of both $so(5)$ and $so(5,1)$ as Lie algebras of 
twisted derivations (having however always in the back of our mind their origin 
as deformed Hopf algebras).

The eight roots of the Lie algebra $so(5)$ are two-component vectors 
$r=(r_1,r_2)$ of the form $r=(\pm1,\pm1)$ and $r=(0,\pm1), r=(\pm1,0)$.  There 
are corresponding generators $E_r$ of $so(5)$ together with two mutually 
commuting generators $H_1,H_2$ of the Cartan subalgebra. The Lie brackets are 
\bea\label{lie-so5} && [H_1,H_2] = 0, \quad [H_j,E_r] = r_j E_r , \nn \\ && 
[E_{-r},E_{r}] = r_1 H_1 + r_2 H_2, \quad [E_{r},E_{r'}] = N_{r,r'} E_{r+r'}, 
\eea with $N_{r,r'}=0$ if $r+r'$ is not a root.

In order to give the action of $so(5)$ on $\S^4$, for convenience we introduce 
``partial derivatives'', $\partial_\mu$ and $\partial_\mu^*$ with the usual 
action on the generators of the algebra $\A(\S^4)$ {i.e},
 $\partial_\mu(z_\nu)=\delta_{\mu\nu}$, $\partial_\mu(z_\nu^*)=0$, and
 $\partial_\mu^*(z_\nu^*)=\delta_{\mu\nu}$, $\partial_\mu^*(z_\nu)=0$. With these 
we construct operators on $\A(\S^4)$, \be\label{act4} \begin{aligned} H_1 &= z_1 
\partial_1 - z_1^* \partial_1^* , \\ E_{+1,+1}&= z_2 \partial_1^* - z_1 
\partial_2^*,\\ E_{+1,0} &= \tfrac{1}{\sqrt{2}} (2 z_0 \partial_1^* - z_1 
\partial_0), \end{aligned} \qquad \begin{aligned} H_2 &= z_2 \partial_2 - z_2^* 
\partial_2^* \\ E_{+1,-1} &= z_2^* \partial_1^* - z_1 \partial_2 \, , \\ E_{0,+1} 
&= \tfrac{1}{\sqrt{2}} (2 z_0 \partial_2^* - z_2 \partial_0) \, , \end{aligned} 
\ee and $E_{-r}=(E_{r})^*$, with the obvious meaning of the adjoint.  
\begin{prop}\label{prop:actso5} The operators in \eqref{act4} give a well defined 
action of $so(5)$ on the algebra $\A(\S^4)$ provided one extends them to the 
whole of $\A(\S^4)$ as ``twisted derivations'' via the rules, \bea\label{tder} && 
E_r( a b ) = E_r(a) \lambda^{{\frac{1}{2}}({-r_1 H_2+r_2 H_1})}(b) + 
\lambda^{{\frac{1}{2}}({r_1 H_2-r_2 H_1})}(a) E_r(b) , \nn\\ && H_j( a b ) = 
H_j(a) b + a H_j(b) , \eea for any two elements $a,b\in \A(\S^4)$; here 
$\lambda=e^{2\pi\ii \theta}$ is the deformation parameter. \end{prop} 
\begin{proof} With these twisted rules, one explicitly checks compatibility of 
the action \eqref{act4} with the commutation relations \eqref{s4t} of $\A(\S^4)$. 
\end{proof} \begin{rem} The operators $\lambda^{\pm \half r_i H_j}$ in 
\eqref{tder} are understood as exponentials of diagonal matrices: on the 
generators $z_\mu,z_\mu^*$ of $\A(\S^4)$, the operators $H_1$ and $H_2$ can be 
written as finite dimensional matrices. A comparison with 
equation~\eqref{eq:act-S4} shows that $H_1$ and $H_2$ in~\eqref{act4} are the 
infinitesimal generators of the action of $\bT^2$ on $\S^4$.  \end{rem}

We can write the twisted action of $so(5)$ on $\A(\S^4)$ by using the 
quantization map $L_\theta$ introduced in Section~\ref{section:toric-ncm}. For 
$L_\theta(a) \in \A(\S^4)$ and $t \in so(5)$ a twisted action is defined by \be 
\label{eq:def-twisted} T \cdot L_\theta (a) = L_\theta(t\cdot a) \ee where $T$ is 
the ``quantization'' of $t$ 
and $t \cdot a$ is the classical action of $so(5)$ on 
$\A(S^4)$ (a better but heavier notation for the action $T~\cdot$ would be 
$t~\cdot_\theta $). One checks that both of these definitions of the twisted action 
coincide. The latter definition allows one to define an action of $so(5)$ on 
$\Cinf(\S^4)$ by allowing $a$ to be in $\Cinf(S^4)$ in 
equation~\eqref{eq:def-twisted}. Furthermore, as operators on the Hilbert space 
$\ch$ of spinors, one could identify $\lambda^{{\frac{1}{2}}({r_1 H_2-r_2 
H_1})}=U(\frac{1}{2} r\cdot \theta)$, with $r=(r_1,r_2)$, $\theta$ the 
antisymmetric two by two matrix with $\theta_{12}=-\theta_{21}=\theta$ and $U(s)$ 
is the representation of $\IT^2$ on $\ch$ as in Sect~\ref{subsection:def-torus}.

The twisted action of the Lie algebra $so(5)$ on $\A(\S^4)$ is extended to the 
differential calculus $(\Omega(\S^4),\dd)$ by requiring it to commute with the 
exterior derivative, $$ T \cdot \dd \omega := \dd (T \cdot \omega) . $$ for $T\in 
so(5), ~ \omega \in \Omega(\S^4)$. Then, we need to use the rule \eqref{tder} on 
a generic form. For instance, on 1-forms we have, \begin{align}\label{onforms} & 
E_r( \sum_k a_k \dd b_k ) = \sum_k \Big( E_r(a_k) \dd \big( 
\lambda^{{\frac{1}{2}}({-r_1 H_2+r_2 H_1})}(b_k) \big) + 
\lambda^{{\frac{1}{2}}({r_1 H_2-r_2 H_1})}(a_k) \dd \big( E_r(b_k) \big) \Big), 
\nn \\ & H_j( \sum_k a_k \dd b_k ) = \sum_k \Big(H_j(a_k) \dd b_k + a_k \dd 
\big(H_j (b_k) \big) \Big)\, . \end{align}

The representation of $so(5)$ on $\S^4$ given in \eqref{act4} is the fundamental 
vector representation. When lifted to $\Sk^7$ one gets the fundamental spinor 
representation: as we see from the
 quadratic relations among corresponding generators, as given in 
\eqref{subalgebra}, the lifting amounts to take the ``square root'' 
representation.
 The action of $so(5)$ on $\A(\Sk^7)$ is constructed by requiring twisted 
derivation properties via the rule \eqref{tder} -- when acting now on any two 
elements $a,b\in \A(\S^7)$ -- so as to reduce to the action \eqref{act4} on 
$\A(\S^4)$ when using the defining
 quadratic relations \eqref{subalgebra}. The resulting action on $\A(\Sk^7)$ can 
be given as the action of matrices $\Gamma$'s on the $\psi$'s, 
\begin{align}\label{act7} \psi_a \mapsto \sum_b \Gamma_{ab} \psi_b, \qquad 
\psi^*_a \mapsto \sum_b \tilde\Gamma_{ab} \psi^*_b \, , \end{align} with the 
matrices $\Gamma = \{H_j, E_r\}$ given explicitly by, \be\label{tgamma} 
\begin{aligned} &H_1 = \half\left( \begin{smallmatrix} 1 & & & \\
 & -1 & & \\
 & & -1 & \\
 & & & 1 \end{smallmatrix} \right), \\ &E_{+1,+1} =\begin{pmatrix} 0 & 0 \\ 0 & 
\begin{smallmatrix} 0 & -1 \\ 0 & 0 \end{smallmatrix} \\ \end{pmatrix}, \\ 
&E_{+1,0} = \tfrac{1}{\sqrt{2}} \begin{pmatrix} 0 & \begin{smallmatrix} 0 & 0\\ 0 
& -1 \end{smallmatrix} \\ \begin{smallmatrix} \mu & 0 \\ 0 & 0 \end{smallmatrix} 
& 0 \\ \end{pmatrix}, \end{aligned} \qquad \begin{aligned} &H_2 = \half\left( 
\begin{smallmatrix} -1 & & & \\
 & 1 & & \\
 & & -1 & \\
 & & & 1 \end{smallmatrix} \right), \\ &E_{+1,-1} =\begin{pmatrix} 
\begin{smallmatrix} 0 & 0 \\ -\mu & 0 \end{smallmatrix} & 0 \\ 0 & 0 \\ 
\end{pmatrix}, \\ &E_{0,+1} = \tfrac{1}{\sqrt{2}} \begin{pmatrix} 0 & 
\begin{smallmatrix} 0 & \bar{\mu} \\ 0 & 0 \end{smallmatrix} \\ 
\begin{smallmatrix} 0 & 1 \\ 0 & 0 \end{smallmatrix} & 0 \\ \end{pmatrix}, 
\end{aligned} \ee and $\tilde\Gamma:= \sigma \Gamma \sigma^{-1}$ with \be \sigma 
:= \begin{pmatrix} \begin{smallmatrix} 0 & -1 \\ 1 & 0 \end{smallmatrix} & 0 \\ 0 
& \begin{smallmatrix} 0 & -1 \\ 1 & 0 \end{smallmatrix} \end{pmatrix}.\ee 
Furthermore, $E_{-r}=(E_{r})^*$. With the twisted rules \eqref{tder} for the 
action on products of any two elements $a,b\in \A(\S^7)$, one checks 
compatibility of the above action with the commutation relations \eqref{s7t} of 
$\A(\Sk^7)$. Again, the operators $\lambda^{\pm \half r_i H_j}$ in \eqref{tder} 
are exponentials of diagonal matrices $H_1$ and $H_2$ given in the representation 
\eqref{tgamma} and as above, one could think of $\lambda^{{\frac{1}{2}}({r_1 
H_2-r_2 H_1})}$ as the operator $U(\frac{1}{2} r\cdot \theta)$. \begin{rem} 
Compare the form of the matrices $H_1$ and $H_2$ in the representation 
\eqref{tgamma} above with the lifted action $\tilde \sigma$ of $\tilde \bT^2$ on 
$\Sk^7$ as defined in \eqref{eq:lift-S7}. One checks that \[ \tilde\sigma_s = e^{ 
\pi i \left((s_1+s_2) H_1 + (-s_1+s_2) H_2\right)} , \] when acting on the spinor 
$(\psi_1, \cdots, \psi_4)$. \end{rem}

Notice that $\tilde\Gamma=-\Gamma^t$ at $\theta=0$. There is a beautiful 
correspondence between the matrices in the representation \eqref{tgamma} and the 
twisted Dirac matrices introduced in \eqref{eq:dirac}, \be \begin{aligned} 
&\tfrac{1}{4}[\gamma_1^*, \gamma_1]=2 H_1\\ &\tfrac{1}{4}[\gamma_1, 
\gamma_2]=(\mu+\bar\mu) E_{+1,+1}\\ &\tfrac{1}{4}[\gamma_1, \gamma_0]= \sqrt{2} 
E_{+1,0} \end{aligned} \qquad \begin{aligned} &\tfrac{1}{4}[\gamma_2^*, 
\gamma_2]=2 H_2\\ &\tfrac{1}{4}[\gamma_1, \gamma_2^*]=(\mu+\bar\mu) E_{+1,-1}\\ 
&\tfrac{1}{4}[\gamma_2, \gamma_0]=\sqrt{2}\bar\mu E_{0,+1} . \end{aligned} \ee 
\begin{rem} It is straightforward to check that the twisted Dirac matrices 
satisfy the following relations under conjugation by $\sigma$: \be 
\label{eq:gamma-gammat} (\sigma\gamma_0\sigma^{-1})^t = \gamma_0;\qquad 
(\sigma\gamma_1\sigma^{-1})^t = \gamma_1 \lambda^{H_2}; \qquad 
(\sigma\gamma_2\sigma^{-1})^t = \gamma_2 \lambda^{H_1}. \ee \end{rem}

As for $\S^4$, the twisted action of $so(5)$ on $\A(\S^7)$ is straightforwardly 
extended to the differential calculus $(\Omega(\S^7),\dd)$. Furthermore, due to 
the form of $\tilde\Gamma$ and the property $\Psi_{a2} = \sigma_{ab} \psi^*_b$ 
for the second column of the matrix $\Psi$ in \eqref{Psi}, we have also that 
$so(5)$ acts on $\Psi$ by left matrix multiplication by $\Gamma$, and by right 
matrix multiplication on $\Psi^*$ by the matrix transpose $\tilde \Gamma^t$ as 
follows \be \Psi_{ai} \mapsto \sum_b \Gamma_{ab} \Psi_{bi}, \qquad \Psi^*_{ia} 
\mapsto \sum_a \Psi^*_{ib} \tilde\Gamma_{ab}. \ee All this is used in the 
following \begin{prop} The instanton gauge potential $\omega$ is invariant under 
the action of $so(5)$. \end{prop} \begin{proof} From the above observations, the 
gauge potential transforms as, $$ \omega = \Psi^* \dd \Psi \mapsto \Psi^* \big( 
\tilde\Gamma^t \lambda^{-r_1 H_2} + \lambda^{r_2 H_1} \Gamma \big) \dd \Psi. $$ 
where $\lambda^{-r_i H_j}$ is understood in its representation \eqref{tgamma} on 
$\A(\S^7)$. Direct computation for $\Gamma=\{H_j, E_r\}$ shows that 
$\tilde\Gamma^t \lambda^{-r_1 H_2} + \lambda^{r_2 H_1} \Gamma=0$, which finishes 
the proof. \end{proof}

\subsection{Twisted conformal transformations}\label{subsect:infconf}

In order to have new instantonic configurations we need to use conformal 
transformations. The conformal Lie algebra $so(5,1)$ consists of the generators 
of $so(5)$ together with the dilation and the so-called special conformal 
transformations. On $\bR^4$ with coordinates $\{x_\mu, \mu=1,\ldots,4 \}$ they 
are given by the operators $H_0=\sum_\mu x_\mu \partial/\partial x_\mu$ and 
$G_\mu = 2 x_\mu \sum_\nu x_\nu \partial/\partial x_\nu - \sum_\nu x_\nu^2 
(\partial/\partial x_\nu)$, respectively \cite{Lie70}.

Since we are deforming with respect to the two-torus coming from the Cartan 
subalgebra $H_1,H_2$ of $so(5)$, we write the Lie algebra $so(5,1)$ as $so(5)$ 
together with five extra generators, $H_0, G_r$, the latter labeled by the 
corresponding roots with respect to $H_1$ and $H_2$, that is $r=(\pm 1,0),(0,\pm 
1)$. Besides the Lie brackets of $so(5)$ given in \eqref{lie-so5}, we have 
\begin{align}\label{lie-confa} \begin{aligned} &[H_0,H_i]=0, \\ & 
[H_0,G_r]=\sqrt{2} E_r, \end{aligned} \qquad \begin{aligned} &[H_j, G_r] = r_j 
G_r, \\ & [H_0,E_r]=(\sqrt{2})^{-1} G_r, \end{aligned} \end{align} whenever 
$r=(\pm 1,0),(0,\pm 1)$, and \begin{align}\label{lie-confb} \begin{aligned} 
&[G_{-r},G_r]=2 r_1 H_1 + 2 r_2 H_2, \\ & [E_r,G_{r'}]=\tilde N_{r,r'} G_{r+r'}, 
\end{aligned} \qquad \begin{aligned} &[G_r,G_{r'}] =N_{r,r'} E_{r+r'}, \\ 
&[E_{-r},G_r]=\sqrt{2} H_0, \end{aligned} \end{align} with, as before, 
$N_{r,r'}=0$ if $r+r'$ is not a root of $so(5)$ and $\tilde N_{r,r'}=0$ if $ r+r' 
$ does not belong to $\{ (\pm 1,0),(0,\pm 1) \}$.

The action of $so(5,1)$ on $\A(\S^4)$ is given by the operators \eqref{act4} 
together with \begin{align}\label{confact4} &H_0 = \partial_0 - z_0 (z_0 
\partial_0 + z_1 \partial_1 + z_1^* \partial_1^* +
 z_2 \partial_2 + z_2^* \partial_2^*), \nn \\ &G_{1,0}= 2 \partial_1^* - z_1 (z_0 
\partial_0 + z_1 \partial_1 + z_1^* \partial_1^* + \bar\lambda z_2 \partial_2 + 
\lambda z_2^* \partial_2^*),\\ &G_{0,1}= 2 \partial_2^* - z_2 (z_0 \partial_0 + 
z_1 \partial_1 + z_1^* \partial_1^* + z_2 \partial_2 + z_2^* \partial_2^*), \nn 
\end{align} and $G_{-r} = (G_r)^*$. Note that the introduction of the extra 
$\lambda$'s in $G_{1,0}$ (and $G_{-1,0}$) are necessary for the Lie algebra 
structure of $so(5,1)$ -- as dictated by the Lie brackets in \eqref{lie-confa} 
and \eqref{lie-confb} -- to be preserved, that is in order to have a Lie algebra 
representation. Since the operators $H_0$ and $G_r$ are quadratic in the $z$'s, 
one has to be careful when deriving the above Lie brackets and use the twisted 
rules \eqref{tder}. For instance, on the generator $z_2$, we have \begin{align*} 
[E_{-1,-1},G_{1,0}](z_2)&=E_{-1,-1}(-\bar\lambda z_1 z_2)+G_{1,0}(z_1^*)\\ 
&=-\bar\lambda (E_{-1,-1}(z_1) \lambda^{H_2} (z_2)+ \lambda^{H_1}(z_1) 
E_{-1,-1}(z_2))+ G_{1,0}(z_1^*) \\ &=-z_2^* z_2 + z_1 z_1^* +2 - z_1 z_1^* = 
G_{0,-1}(z_2) \end{align*} As for Proposition~\ref{prop:actso5}, a direct 
computation establishes the following, \begin{prop} The operators in \eqref{act4} 
and \eqref{confact4} give a well defined action of $so(5,1)$ on the algebra 
$\A(\S^4)$ provided one extends them to the whole of $\A(\S^4)$ as twisted 
derivations via the rules \eqref{tder} together with \bea && G_r( a b ) = G_r(a) 
\lambda^{{\frac{1}{2}}({-r_1 H_2+r_2 H_1})}(b) + \lambda^{{\frac{1}{2}}({r_1 
H_2-r_2 H_1})}(a) G_r(b) , \nn\\ && H_0( a b ) = H_0(a) b + a H_0(b) , \eea for 
any two elements $a,b\in \A(\S^4)$. \end{prop}

Equivalently, $so(5,1)$ could be defined to act on $\A(\S^4)$ by 
\be\label{action:hopf-so51} T \cdot L_\theta(a) = L_\theta(t\cdot a) , \ee for 
$T$ the operator deforming $t \in so(5,1)$ and $L_\theta(a) \in \A(\S^4)$ 
deforming $a \in \A(S^4)$. Again, equation~\eqref{action:hopf-so51} makes sense 
for $a \in \Cinf(S^4)$, which defines an action of $so(5,1)$ on $\Cinf(\S^4)$. As 
before, the action on the differential calculus $(\Omega(\S^4),\dd)$ is obtained 
by requiring commutation with the exterior derivative: $ T \cdot \dd \omega = \dd 
(T \cdot \omega)$, for $T\in so(5,1)$ and $\omega \in \Omega(\S^4)$. On products 
we shall have formul{\ae} like the one in \eqref{onforms}, \begin{align} & G_r( 
\sum_k a_k \dd b_k ) = \sum_k \Big( G_r(a_k) \dd \big( 
\lambda^{{\frac{1}{2}}({-r_1 H_2+r_2 H_1})}(b_k) \big) + 
\lambda^{{\frac{1}{2}}({r_1 H_2-r_2 H_1})}(a_k) \dd \big( G_r(b_k) \big) \Big), 
\nn \\ & H_0( \sum_k a_k \dd b_k ) = \sum_k \Big(H_0(a_k) \dd b_k + a_k \dd 
\big(H_0 (b_k) \big) \Big)\, . \end{align} What we are dealing with are 
``infinitesimal'' twisted conformal transformations: \begin{lma} 
\label{lma:so51-hodge} The Hodge $\ast_\theta$-structure of $\Omega(\S^4)$ is 
invariant for the twisted action of $so(5,1)$, $$ T\cdot (\ast_\theta\omega)= 
\ast_\theta (T \cdot \omega) , $$ \end{lma} \begin{proof} Recall that 
$T(L_\theta(a))=L_\theta(t\cdot a)$ for $a \in \A(S^4)$ and $T$ the 
``quantization'' of $t \in so(5,1)$. Then, since $so(5,1)$ leaves the Hodge 
$\ast$-structure of $\Omega(S^4)$ invariant and the differential $\dd$ commutes 
with the action of $so(5,1)$, if follows that the latter algebra leaves the Hodge 
$\ast_\theta$-structure of $\Omega(\S^4)$ invariant as well. \end{proof} Again, 
the action of $so(5,1)$ on $\S^4$ can be lifted to an action on $\Sk^7$. And the 
latter action can be written as in \eqref{act7} in terms of matrices $\Gamma$'s 
acting on the $\psi$'s, \begin{align} \psi_a \mapsto \sum_b \Gamma_{ab} \psi_b, 
\qquad \psi^*_a \mapsto \sum_b \tilde\Gamma_{ab} \psi^*_b , \end{align} where in 
addition to \eqref{tgamma} we have also the matrices $\Gamma=\{H_0,G_r\}$, given 
explicitly by \begin{align} H_0 &= \half (-z_0 \I_4 + \gamma_0), \nn\\ 
G_{1,0}&=\half (-z_1 \lambda^{-H_2} + \gamma_1), \qquad G_{0,1}=\half (-z_2 + 
\lambda^{-H_1} \gamma_2), \end{align} with $G_{-r}= (G_r)^*$ and $\tilde\Gamma = 
\sigma \Gamma \sigma^{-1}$. Notice the reappearance of the twisted Dirac matrices 
$\gamma_\mu, \gamma_\mu^*$ of \eqref{eq:dirac} in the above expressions. In the 
above expressions, the operators $\lambda^{-H_j}$ are $4 \times 4$ matrices 
obtained from the spin representation \eqref{tgamma} of $H_1$ and $H_2$ and 
given explicitily by 
\begin{align}\label{expspin} 
\lambda^{-H_1} = \left( \begin{smallmatrix} \bar\mu & & & \\
 & \mu & & \\
 & & \mu & \\
 & & & \bar\mu \end{smallmatrix} \right), \qquad \lambda^{-H_2} = \left( 
\begin{smallmatrix} \mu & & & \\
 & \bar\mu & & \\
 & & \mu & \\
 & & & \bar\mu \end{smallmatrix} \right), \qquad \mu=\sqrt{\lambda}. 
\end{align} 
As for $so(5)$, the action of $so(5,1)$ on the matrix $\Psi$ is found to be by 
left matrix multiplication by $\Gamma$ and on $\Psi^*$ by $\tilde \Gamma$, \be 
\Psi_{ai} \mapsto \sum_b \Gamma_{ab} \Psi_{bi}, \qquad \Psi^*_{ia} \mapsto \sum_a 
\tilde\Gamma_{ab} \Psi^*_{ib}. \ee Here we have to be careful with the ordering 
between $\tilde\Gamma$ and $\Psi^*$ in the second term since the $\tilde\Gamma$'s 
involve the (not-central) $z$'s. 
There are the following useful commutation 
relations between the $z_\mu$'s and $\Psi$: \be \label{eq:z-Psi} \begin{aligned} 
&z_1 \Psi_{ai} = (\lambda^{-H_2})_{ab} \Psi_{bi} z_1, \qquad z_2 \Psi_{ai} = 
(\lambda^{-H_1})_{ab} \Psi_{bi} z_2, \\ &z_1 \Psi^*_{ia} = \Psi^*_{ib} 
(\lambda^{-H_2})_{ba} z_1, \qquad z_2 \Psi^*_{ia} = \Psi^*_{ib} 
(\lambda^{-H_1})_{ba} z_2, 
\end{aligned} \ee
with $\lambda^{-H_j}$ understood as the explicit matrices \eqref{expspin}.

\begin{prop} \label{prop:instantons} The instanton gauge potential $\omega=\Psi^* 
\dd \Psi$ transforms under the action of the extra elements of $so(5,1)$ as 
$\omega \mapsto \omega+ \delta_i \omega$, where \begin{align*} &\delta_0 \omega 
:= H_0(\omega)=-z_0 \omega - \half d z_0 \I_2 + \Psi^* ~\gamma_0 ~\dd\Psi,\\ 
&\delta_1 \omega :=G_{+1,0}(\omega)= -z_1 \omega - \half d z_1 ~\I_2 + \Psi^* 
~\gamma_1~ \dd\Psi,\\ &\delta_2 \omega :=G_{0,+1}(\omega)= -z_2 \omega - \half d 
z_2 ~ \I_2 + \Psi^* ~\gamma_2~ \dd\Psi,\\ &\delta_3 \omega :=G_{-1,0}(\omega)= - 
\omega \bar z_1-\half d \bar z_1 ~\I_2+\Psi^*~\gamma_1^* ~\dd\Psi,\\ &\delta_4 
\omega :=G_{0,-1}(\omega)= -\omega \bar z_2- \half d \bar z_2 ~ 
\I_2+\Psi^*~\gamma_2^*~ \dd\Psi, \end{align*} with $\gamma_\mu, \gamma_\mu^*$ the 
twisted $4\times 4$ Dirac matrices defined in \eqref{eq:dirac}. \end{prop} 
\begin{proof} The action of $H_0$ on the instanton gauge potential $\omega=\Psi^* 
\dd \Psi$ takes the form \begin{align*} H_0(\omega) &=H_0(\Psi^*) \dd \Psi + 
\Psi^* \dd(H_0(\Psi))= \Psi^* (- z_0 \I_4 + \gamma_0 )\dd \Psi - \half d z_0 
\Psi^* \Psi, \end{align*} since $z_0$ is central. Direct computation results in 
the above expression for $\delta_0 \omega $. Instead, the twisted action of $G_r$ 
on $\omega$ takes the form, $$ G_r: \omega_{ij} \mapsto \sum_{a,b,c} 
\tilde\Gamma_{ab} \Psi^*_{ib} (\lambda^{-r_1 H_2})_{ac} \dd \Psi_{cj} 
+(\lambda^{r_2 H_1})_{ab} \Psi^*_{ib} \Gamma_{ac} \dd \Psi_{cj} + (\lambda^{r_2 
H_1})_{ab} \Psi^*_{ib} (\dd \Gamma_{ac}) \Psi_{cj} , $$ where we used the fact 
that $\tilde H_j= \sigma H_j \sigma^{-1} = - H_j$. Let us consider the case 
$r=(+1,0)$. Firstly, note that the complex numbers $(\lambda^{-H_2})_{ac}$ 
commute with $\Psi^*_{ib}$ so that from the definition of $\Gamma$ and 
$\tilde\Gamma$, and using \eqref{eq:z-Psi}, we obtain for the first two terms, 
\begin{align*} -z_1 (\Psi^* \dd \Psi)_{ij} + \half \Psi^*_{ib} 
(\sigma\gamma_1\sigma^{-1})_{cb} (\lambda^{-H_2})_{cd} \dd \Psi_{dj} +\half 
\Psi^*_{ib} (\gamma_1)_{bc} \dd \Psi_{cj}. \end{align*} The first term forms the 
matrix $-z_1 \omega$ whereas the second two terms combine to give $\Psi^* 
\gamma_1 \dd \Psi$, due to relation \eqref{eq:gamma-gammat}. Finally, using 
equation~\eqref{eq:z-Psi} the term $ \Psi^*_{ib} (\dd \Gamma_{ac}) \Psi_{cj}$ 
reduces to $-\half\dd z_1 \Psi^*_{ib} \Psi_{bj}=-\half \dd z_1 \I_2$. The 
formulae for $r=(-1,0)$ and $r=(0, \pm 1)$ are established in likewise manner. 
\end{proof} \begin{rem} At first sight, the infinitesimal gauge potentials 
$\delta_j \omega$ given above do not seem to be $su(2)$-gauge potentials, in that 
they do not satisfy $\bar{(\delta_j \omega)_{kl}}= (\delta_j \omega)_{lk}$ and 
$\sum_k (\delta_j \omega)_{kk}=0$. This is only due to the fact that the 
generators $G_r$ and $H_0$ are the deformed analogues of the generators of the 
complexified Lie algebra $so(5,1) \ot_\bR \C$. One recovers $su(2)$-gauge 
potentials by acting with the real generators $\half(G_r + G_r^*)$, $\frac{1}{2i} 
(G_r - F_r^*)$ and $H_0$. The resulting gauge potentials, $\delta_0\omega$, 
$\half(\delta_1 \omega + \delta_3 \omega)$, $\frac{1}{2i}(\delta_1 \omega - 
\delta_3 \omega)$, $\half(\delta_2 \omega + \delta_4 \omega)$ and 
$\frac{1}{2i}(\delta_2 \omega - \delta_4 \omega)$, are traceless skew-hermitian 
matrices with entries in $\Omega^1(\Sk^7)$. \end{rem} The transformations of the 
gauge potential $\omega$ under the twisted symmetry $so(5,1)$, given in 
Proposition~\ref{prop:instantons}, induce natural transformations of the 
canonical connection $\nabla_0$ in \eqref{cancon} to $\nabla_{t,i}:=\nabla_0 + t 
\delta_i \omega + \cO(t^2)$. We shall presently see that these new connections 
are (infinitesimal) instantons, {i.e.} their curvatures are self-dual. In fact, 
this also follows from Lemma~\ref{lma:so51-hodge} which states that $so(5,1)$ 
acts by conformal transformation therefore leaving the self-duality equation 
$\ast_\theta F_0 = F_0$ for the basic instanton $\nabla_0$ invariant.

We start by writing $\nabla_{t,i}$ in terms of the canonical connection on 
$\E\isom p \big(\A(\S^4)\big)^4$. Using the isomorphism, described in 
Example~\ref{ex:instanton-bundle}, between this module and the module of 
equivariant maps $\A(\Sk^7)\boxtimes_\rho \C^2$, we find that $\nabla_{t,i}=p \dd 
+ t \delta_i \alpha + \cO(t^2)$ with \be \begin{aligned} \delta_0 \alpha = p 
\gamma_0 (\dd p) p - \half \Psi \dd z_0 \Psi^*, \\ \delta_1 \alpha = p \gamma_1 
(\dd p) p - \half \Psi \dd z_1 \Psi^*,\\ \delta_2 \alpha = p \gamma_2 (\dd p) p - 
\half \Psi \dd z_2 \Psi^*, \end{aligned} \qquad \begin{aligned} ~ \\ \delta_3 
\alpha = p \gamma_1^* (\dd p) p - \half \Psi \dd z_1^* \Psi^*,\\ \delta_4 \alpha 
= p \gamma_2^* (\dd p) p - \half \Psi \dd z_2^* \Psi^*, \end{aligned} \ee The 
$\delta_i \alpha$'s are $4 \times 4$ matrices with entries in the one-forms 
$\Omega^1(\S^4)$ and satisfying conditions $p \delta_i \alpha = \delta_i \alpha p 
= p \delta_i \alpha p = \delta_i \alpha$, as expected from the general theory on 
connections on modules in Section~\ref{section:connections}. Indeed, using 
relations \eqref{eq:z-Psi} one can move the $\dd z_\mu$'s to the left of $\Psi$ 
at the cost of some $\mu$'s, so getting expression like $\dd z_i ~p \in 
M_4(\Omega^1(\S^4))$.

From equation~\eqref{ucurv}, the curvature $F_{t,i}$ of the connection 
$\nabla_{t,i}$ is given by \be F_{t,i}=F_0 + t p \dd (\delta_i \alpha) + 
\cO(t^2). \ee In order to check self-duality (modulo $t^2$) of this curvature, we 
will express it in terms of the projection $p$ and consider $F_{t,i}$ as a 
two-form valued endomorphism on the module $\E\isom p \big(\A(\S^4)\big)^4$. \begin{prop} 
\label{prop:curv-transf} The curvatures $F_{t,i}$ of the connections 
$\nabla_{t,i}$, $i=0,\ldots,4$, are given by $F_{t,i} =F_0 + t \delta_i F+ 
\cO(t^2)$, where $F_0=p \dd p \dd p$ and $\delta_i F$ are the following $4\times 
4$-matrices of 2-forms: \be \begin{aligned} \delta_0 F&=-2 z_0 F_0, \\ \delta_1 
F&=-2 z_1 \lambda^{H_2} F_0, \\ \delta_2 F&=-2 z_2 \lambda^{H_1} F_0; 
\end{aligned} \qquad \begin{aligned} ~ \\ \delta_3 F=-2 z_1^*\lambda^{-H_2} F_0, 
\\ \delta_4 F=-2 z_2^*\lambda^{-H_1} F_0. \end{aligned} \ee \end{prop} 
\begin{proof} A small computation yields for $\delta_i F=p \dd (\delta_i 
\alpha)$, thought of as an $\Omega^2(\S^4)$-valued endomorphism on $\E$ the 
expression, $ \delta_i F = p (\dd p) \gamma_i (\dd p) p - p \gamma_i (\dd p) (\dd 
p) p, $ with the notation $\gamma_3=\gamma_1^*$ and $\gamma_4=\gamma_2^*$, and 
using $p (\dd p) p=0$. Then, the crucial property $p (\dd p\gamma_i +\gamma_i \dd 
p) (\dd p) p = 0$ all $i=0,\ldots,4$ yields $\delta_i F = -2 p \gamma_i \dd p \dd 
p p$. This is expressed as $\delta_i F = -2 p \gamma_i p dp dp$ by using $\dd p = 
(\dd p) p + p \dd p $. Finally, $p \gamma_i p=\Psi (\Psi^* \gamma_i \Psi) 
\Psi^*$, so that the result follows from the definition of the $z$'s in terms of 
the Dirac matrices given in equation~\eqref{zg}, together with the commutation 
relations between them and the matrix $\Psi$ in equation~\eqref{eq:z-Psi} 
\end{proof} \begin{prop} The connections $\nabla_{t,i}$ are (infinitesimal) 
instantons, {i.e.} \be \ast_\theta F_{t,i}=F_{t,i} \qquad \mathrm{mod} \, t^2 . 
\ee Moreover, the connections $\nabla_{t,i}$ are not gauge equivalent to 
$\nabla_0$.\end{prop} \begin{proof} The first point follows directly from the 
above expressions for $\delta_i F$ and the self-duality of $F_0$.  To establish 
the gauge inequivalence of the connections $\nabla_{t,i}$ with $\nabla_0$, we 
recall that an infinitesimal gauge transformation is given by $\nabla_0 \mapsto 
\nabla_0 + t [\nabla_0,X]$ for $X \in \Gamma^\infty(\ad(\Sk^7))$. We need to show 
that $\delta_i \omega$ is orthogonal to $[\nabla_0,X]$ for any such $X$, {i.e.} 
$$ ([\nabla_0,X],\delta_i \omega )_2= 0, $$ with the natural inner product on 
$\Omega^1(\ad(\Sk^7)):=\Omega^1(\S^4)\ot_{\Cinf(\S^4)} 
\Gamma^\infty(\ad(\Sk^7))$. From Remark~\ref{rem:conn-form-central-general}, it 
follows that $$(\nabla_0^{(2)}(X),\delta_i \omega)_2 = (X, 
\left(\nabla_0^{(2)}\right)^* (\delta_i \omega) )_2,$$ which then should vanish 
for all $X$. From equation~\eqref{action:hopf-so51}, we see that $\delta_i \omega 
= T_i(\omega)$ coincides with $L_\theta (t_i \cdot \omega^\class)$ with $t_i$ and 
$\omega^\class$ the classical counterparts of $T_i$ and $\omega$, respectively. 
In the undeformed case, the infinitesimal gauge potentials generated by acting 
with elements in $so(5,1)-so(5)$ on the basic instanton gauge potential 
$\omega^\class$ satisfy $(\nabla_0^{(2)} )^* (\delta_i \omega^\class)=0$ as shown 
in \cite{AHS78}. The result then follows from the observation that 
$\nabla_0^{(2)}$ commutes with the quantization map $L_\theta$ (cf. 
Remark~\ref{rem:conn-form-central-general}). \end{proof}

\subsection{Local expressions} \label{subsect:localexpr} In this section, we 
obtain ``local expressions'' for the instantons on $\S^4$ constructed in the 
previous section; that is we map them to a noncommutative $\R^4$ obtained by 
``removing a point'' from $\S^4$. The algebra $\A(\R^4)$ of polynomial functions 
on the 4-plane $\R^4$ is defined to be the $*$-algebra generated by elements 
$\zeta_1, \zeta_2$ satisfying \be \zeta_1 \zeta_2 = \lambda \zeta_2 \zeta_1; 
\qquad \zeta_1 \zeta_2^* = \bar\lambda \zeta_2^* \zeta_1. \ee with 
$\lambda=e^{2\pi i \theta}$ as above. At $\theta=0$ one recovers the $*$-algebra 
$\A(\bR^4)$ of polynomial functions on the usual 4-plane $\bR^4$.

 The algebra $\A(\R^4)$ can also be defined as the vector space $\A(\bR^4)$
 equipped with a deformed product $\times_\theta$ as in
 equation~\eqref{eq:star-product}. Indeed, the torus $\bT^2$ acts naturally on
 the two complex coordinates of $\bR^4 \isom \C^2$. This also allows us to define
 the smooth algebra $C_b^\infty(\R^4)$ as the vector space $C_b^\infty(\bR^4)$ of
 bounded smooth functions on $\bR^4$ equipped with a deformed product
 $\times_\theta$. However, for our purposes it suffices to consider the
 polynomial algebra $\A(\R^4)$ with one self-adjoint central generator $\rho$
 added together with relations $\rho^2 (1+|\zeta|^2)=(1+|\zeta|^2) \rho^2 = 1 $
 where $|\zeta|^2:=\zeta_1^* \zeta_1 + \zeta_2^* \zeta_2$ (this enlargement was
 already done in \cite{CD02}). In the following, we will denote the enlarged
 algebra by $\tilde\A(\R^4)$ and will also use the notation 
\be
\rho^2=(1+|\zeta|^2)^{-1}= \frac{1}{1+|\zeta|^2}.
\ee 
Note that $\rho^2$ is an element in $C^\infty_b(\R^4)$.

One defines elements $\tilde z_\mu, \, \mu=0,1,2$ in $\tilde\A(\bR^4)$ by \be 
\label{chart} \tilde z_0 = (1-|\zeta|^2)(1+|\zeta|^2)^{-1}, \qquad \tilde z_j = 2 
\zeta_j (1+|\zeta|^2)^{-1} \quad j=1,2, \ee
 and checks that they satisfy the same relations as in \eqref{s4t} of the
 generators $z_\mu$ of $\A(\S^4)$. The difference is that the classical point
 $z_0=-1, z_j=z_j^*=0$ of $\S^4$ is not in the spectrum of $\tilde z_\mu$. We
 interpret the noncommutative plane $\R^4$ as a ``chart'' of the noncommutative
 4-sphere $\S^4$ and equation~\eqref{chart} as the (inverse) stereographic
 projection from $\S^4$ to $\R^4$. In fact, one can cover $\S^4$ by two such
 charts with domain $\R^4$, and transition functions on $\R^4 \backslash \{0\}$,
 where $\{0\}$ is the classical point $\zeta_j=\zeta_j^* = 0$ of $\R^4$ (cf. 
\cite{CD02} for more details).

A differential calculus $(\Omega(\R^4),\dd)$ on $\R^4$ is obtained from the 
general procedure described in Section~\ref{section:toric-ncm}. Explicitly, 
$\Omega(\R^4)$ is the graded $*$-algebra generated by the elements $\zeta_\mu$ of 
degree $0$ and $d\zeta_\mu$ of degree $1$ with relations, \begin{align} 
\label{rel:diff-R4} &d\zeta_\mu d\zeta_\nu+ \lambda_{\mu\nu} d\zeta_\nu 
d\zeta_\mu =0 , \qquad d\zeta_\mu d\zeta_\nu^* + \lambda_{\nu\mu} d\zeta_\nu^* 
d\zeta_\mu =0, \nn \\ &\zeta_\mu d\zeta_\nu = \lambda_{\mu\nu} d\zeta_\nu 
\zeta_\mu, \qquad \zeta_\mu d\zeta_\nu^* = \lambda_{\nu\mu} d\zeta_\nu^* 
\zeta_\mu , \end{align} and $\lambda_{1 2} = \bar{\lambda}_{2 1} =: 
\lambda=e^{2\pi i \theta}$. There is a unique differential $\dd$ on 
$\Omega(\R^{4})$ for which one has $\dd: \zeta_\mu\mapsto d\zeta_\mu$ and a Hodge 
star operator $\ast_\theta : \Omega^p(\R^4) \to \Omega^{4-p}(\R^4)$, obtained 
from the classical Hodge star operator as before. In terms of the standard 
Riemannian metric on $\bR^4$, on two-forms we have, \begin{align} 
\label{eq:hodge-2forms} \ast_\theta d\zeta_1 d\zeta_2 = -d\zeta_1 d\zeta_2 , 
\qquad \ast_\theta d\zeta_1 d \zeta_1^* = -d\zeta_2 d \zeta_2^* , 
\qquad\ast_\theta d\zeta_1 d \zeta_2^* = d\zeta_1 d \zeta_2^* , \end{align} and 
$\ast_\theta^2 = \id$. These are the same formulae as the ones for the undeformed 
Hodge $\ast$ on $\bR^4$ -- since the metric is not changed in an isospectral 
deformation.

Again, we slightly enlarge the differential calculus $\Omega(\R^4)$ by adding the 
self-adjoint central generator $\rho$. The differential $\dd$ on $\rho$ is 
derived from the Leibniz rule for $\dd$ applied to its defining relation, 
\begin{align*} (\dd \rho^2) (1+|\zeta|^2) + \rho^2 \dd(1+|\zeta|^2)= \dd 
\left(\rho^2 (1+|\zeta|^2) \right) = 0, \end{align*} so that $\rho \dd \rho = 
\half \dd \rho^2 = - \half \rho^4 \dd(1+|\zeta|^2) = - \half \rho^4 \sum_\mu (\zeta_\mu 
d\zeta_\mu^* + \zeta_\mu^* d\zeta_\mu)$. The enlarged differential calculus will 
be denoted by $\tilde\Omega(\R^4)$.

\bigskip The stereographic projection from $S^4$ onto $\bR^4$ is a conformal map 
commuting with the action of $\bT^2$; thus it makes sense to investigate the form 
of the instanton connections on $\S^4$ obtained in 
Proposition~\ref{prop:instantons} on the local chart $\R^4$. As in \cite{Lnd05}, 
we first introduce a ``local section'' of the principal bundle $\Sk^7 \to \S^4$ 
on the local chart of $\S^4$ defined in \eqref{chart}. Let $u=(u_1,u_2)$ be a 
complex spinor of modulus one, $u_1^* u_1 + u_2^*u_2 = 1$, and define \be 
\label{lose} \begin{pmatrix} \psi_1 \\ \psi_2 \end{pmatrix} = \rho 
\begin{pmatrix} u_1 \\ u_2 \end{pmatrix} , \qquad \begin{pmatrix} \psi_3 \\ 
\psi_4 \end{pmatrix} = \rho \begin{pmatrix} \zeta_1^*& \zeta_2^* \\ -\mu \zeta_2 
& \bar{\mu} \zeta_1 \end{pmatrix} \begin{pmatrix} u_1 \\ u_2 \end{pmatrix}. \ee 
\begin{rem} Strictly speaking, the symbols $\psi_a$ here denotes elements in the 
algebra $\A(\Sk^7)$ enlarged by an extra generator which is the inverse of $1+z_0 
= 2(1+\psi_1^* \psi_1 + \psi_2^* \psi_2)$. Intuitively, this corresponds to 
``remove'' the fiber $S^3$ in $\Sk^7$ above the classical point $z_0=-1, 
z_j=z_j^*=0$ of the base space $\S^4$. \end{rem} The commutations rules of the 
$u_j$'s with the $\zeta_k$'s are dictated by those of the $\psi_j$: \be u_1 
\zeta_j = \mu \zeta_j u_1 \, , \quad u_2 \zeta_j = \bar{\mu} \zeta_j u_2 \, , 
\quad j=1,2 . \ee The right action of $\SU(2)$ rotates the vector $u$ while 
mapping to the ``same point'' of $\S^4$, which, using the 
definition~\eqref{subalgebra}, from the choice in \eqref{lose} is found to be \be 
2 (\psi_1 \psi^*_3+ \psi^*_2 \psi_4) = \wt z_1, \, \quad 2(- \psi^*_1 \psi_4 + 
\psi_2\psi^*_3) =\wt{z}_2, \, \quad 2(\psi^*_1 \psi_1 + \psi^*_2 \psi_2) -1 = 
\wt{z}_0, \ee and is in the local chart \eqref{chart}, as expected.

By writing the unit vector $u$ as a matrix, $u=\left(\begin{smallmatrix} u_1 & 
-u_2^* \\ u_2 & u_1^*\end{smallmatrix}\right) \in \SU(2)$, we have \be \Psi= \rho 
\begin{pmatrix} \I_2 & 0 \\ 0 & \cZ \end{pmatrix} \begin{pmatrix} u \\ u 
\end{pmatrix}, \qquad \text{with } \cZ = \left( \begin{smallmatrix} \zeta_1^*& 
\zeta_2^* \\ -\mu \zeta_2 & \bar{\mu} \zeta_1 \end{smallmatrix} \right) . \ee 
Then, direct computation of the gauge potential $\omega= \Psi^* d \Psi$ yields 
\begin{align} u \omega u^* &= \rho^{-1} \dd \rho + \rho^2 \cZ^* \dd \cZ + (\dd u) 
u^* \nn \\ &= \frac{1}{(1+|\zeta|^2)} \begin{pmatrix} \sum_i \zeta_i d \zeta_i^* 
- d \zeta_i \zeta_i^* & 2( \zeta_1 d\zeta_2^* - d \zeta_1 \zeta_2^*) \\ 2(\zeta_2 
d \zeta_1^* - d \zeta_2 \zeta_1^*)  & \sum_i \zeta_i^* d \zeta_i - d \zeta_i^* 
\zeta_i \end{pmatrix} + (\dd u) u^* , \end{align} while its curvature $F_0=\dd 
\omega + \omega^2$ is, \be\label{gcur} u F_0 u^* = \rho^4 \dd\cZ^* \dd\cZ = 
\frac{1}{(1+|\zeta|^2)^2} \begin{pmatrix} d \zeta_1 d \zeta_1^* - d \zeta_2 d 
\zeta_2^* & 2 d \zeta_1 d \zeta_2^* \\ 2 d \zeta_2 d \zeta_1^* & -d \zeta_1 d 
\zeta_1^* - d \zeta_2 d \zeta_2^* \end{pmatrix} . \ee From the 
expressions~\eqref{eq:hodge-2forms} for the Hodge operator on two forms, one 
checks that this curvature is self-dual, $\ast_\theta(u F_0 u^*) = u F_0 u^*$, as 
expected.

\bigskip The explicit local expressions for the transformed -- under 
infinitesimal conformal transformations -- gauge potentials and their curvature 
can be obtained in a similar manner. As an example, let us work out the local 
expression for $\delta_0 \omega$ which is the most transparent one. Given the 
expression for $\delta_0 \omega$ in Proposition~\ref{prop:instantons}, a direct 
computation shows that its local counterpart is \be u \delta_0\omega u^* = - 2 
\rho \dd \rho \, -2 \rho^4 \cZ^* \dd \cZ , \ee giving for the transformed 
curvature, \be u F_{t,0} u^* =F_0 + 2t (1-2\rho^2) F_0 + \cO(t^2). \ee It is 
clear that this rescaled curvature still satisfies the self-duality equation; 
this is also in concordance with Proposition~\ref{prop:curv-transf}, being 
$\tilde z_0 =2\rho^2-1$.

\subsection{Moduli space of instantons}\label{subsect:modulispace} We will 
closely follow the infinitesimal construction of instantons for the undeformed 
case given in \cite{AHS78} . This will eventually result in the computation of 
the dimension of the ``tangent space'' to the moduli space of instantons on 
$\S^4$ by index methods. It will turn out that the five-parameter family of 
instantons constructed in the previous section is indeed the complete set of 
infinitesimal instantons on $\S^4$.

Let us start by considering a family of connections on $\S^4$, \be 
\nabla_t=\nabla_0 + t \alpha \ee where $\alpha \in \Omega^1(\ad(\Sk^7)) \equiv 
\Omega^1(\S^4) \otimes_{\Cinf(\S^4)} \Gamma^\infty(\ad(\Sk^7))$ and 
after Example~\ref{ex:instanton-bundle} we have denoted 
$\Gamma^\infty(\ad(\Sk^7))=\Cinf(\Sk^7)\otimes_\ad su(2)$. 
For $\nabla_t$ to be an instanton, we have to 
impose the self-duality equation $\ast_\theta F_t = F_t$ on the curvature $F_t = 
F_0 + t [\nabla,\alpha] + \cO(t^2)$ of $\nabla_t$. This leads, when 
differentiated with respect to $t$, at $t=0$, to the {\it linearized self-duality 
equation} \be P_- [\nabla_0, \alpha]=0, \ee with $P_-:=\half(1-\ast_\theta)$ the 
projection onto the antiself-dual 2-forms. Here $[\nabla_0, \alpha]$, \be 
[\nabla_0,\alpha]_{ij} = \dd \alpha_{ij} + \omega_{ik} \alpha_{kj} -\alpha_{ik} 
\omega_{kj} , \ee is an element in $\Omega^2(\ad(\Sk^7))$ and has vanishing 
trace, due to the fact that $\omega_{ik} \alpha_{kj} =\alpha_{kj} \omega_{ik}$ 
(cf. equations \eqref{cancon}, \eqref{cangp} and the related discussion).

If the family $\nabla_t$ were obtained from an infinitesimal gauge transformation 
we would have had $\alpha=[\nabla_0,X]$, for some $X \in 
\Gamma^\infty(\ad(\Sk^7))$. Indeed, $[\nabla_0,X]$ is an element in 
$\Omega^1(\ad(\Sk^7))$ and $P_-[\nabla_0,[\nabla_0,X]]=[P_- F_0,X]=0$, since 
$F_0$ is self-dual. Hence, we have defined an element in the first cohomology 
group $H^1$ of the so-called {\it self-dual complex}, \be\label{complex} 0 \to 
\Omega^0(\ad(\Sk^7)) \overset{\dd_0}{\longrightarrow} \Omega^1(\ad(\Sk^7)) 
\overset{\dd_1}{\longrightarrow} \Omega^2_-(\ad(\Sk^7)) \to 0 \ee where 
$\Omega^0(\ad(\Sk^7))=\Gamma^\infty(\ad(\Sk^7))$ and $\dd_0=[\nabla_0,\cdot ~]$, 
$\dd_1:= P_-[\nabla_0,\cdot ~]$. Note that these operators are Fredholm 
operators, so that the cohomology groups of the complex are finite dimensional. 
As usual, the complex can be replaced by a single Fredholm operator \be 
\label{eq:op-fredholm} \dd_0^* + \dd_1 : \Omega^1(\ad(\Sk^7)) \longrightarrow 
\Omega^0(\ad(\Sk^7)) \op \Omega^2_-(\ad(\Sk^7)) , \ee with $\dd_0^*$ the adjoint 
of $\dd_0$ with respect to the inner product \eqref{eq:inner-product-forms}.

Our goal is to compute $h^1=\dim H^1$, the number of ``true'' instantons. This is 
achieved by calculating the alternating sum $h^0-h^1+h^2$ of Betti numbers from 
the index of this Fredholm operator, \be \ind(d_0^* + d_1) = -h^0 + h^1 - h^2, 
\ee while showing that $h^0=h^2=0$.

By definition, $H^0$ consists of the covariant constant elements in 
$\Gamma^\infty(\ad(\Sk^7))$. Since the operator $[\nabla_0,\cdot~]$ commutes with 
the action of $\bT^2$ and coincides with $\nabla_0^{(2)}$ on 
$\Gamma^\infty(\ad(\Sk^7))$ (cf. Remark~\ref{rem:conn-form-central-general}), 
being covariantly constant means that \be [\nabla_0,X]=\nabla_0^{(2)}(X)=0. \ee 
If we write once more $X=L_\theta(X^\class)$ in terms of its classical 
counterpart and use the fact that $\nabla_0^{(2)}$ commutes with $L_\theta$ (cf. 
Remark~\ref{rem:conn-form-central-general}) we find that this condition entails 
\be \nabla_0^{(2)}(L_\theta(X^\class)) = L_\theta( \nabla_0^{(2)}(X^\class))=0. 
\ee For the undeformed case, there are no covariant constant elements in 
$\Gamma^\infty(\ad(S^7))$ for an irreducible self-dual connection on $\E$, thus 
we conclude that $h^0=0$. A completely analogous argument for the kernel of the 
operator $d_1^*$ shows that also $h^2=0$.

\subsection{Dirac operator associated to the complex} The next step consists in 
computing the index of the Fredholm operator $\dd_0^* + \dd_1$ defined in 
\eqref{eq:op-fredholm}. Firstly, this operator can be replaced by a Dirac 
operator on the spinor bundle $\cS$ with coefficients in the adjoint bundle. For 
this, we need the following lemma, which is a straightforward modification of its 
classical analogue \cite{AHS78}. Recall that the $\bZ^2$-grading $\gamma_5$ 
induces a decomposition of the spinor bundle $\cS= \cS^+ \oplus \cS^-$. Note also 
that $\cS^-$ coincides classically with the charge $-1$ anti-instanton bundle. 
Indeed, the Levi-Civita connection -- when lifted to the spinor bundle and 
restricted to negative chirality spinors -- has antiself-dual curvature. 
Similarly, $\cS^+$ coincides with the charge $+1$ instanton bundle. Then 
Remark~\ref{rem:associated-modules} implies that the $\Cinf(S^4)$-modules 
$\Gamma^\infty(S^4,\cS^\pm)$ have a module-basis that is homogeneous under the 
action of $\tilde\bT^2$. We conclude from $\tilde\bT^2$-equivariance that 
$\Gamma^\infty(\S^4, \cS^-)$ is isomorphic to the charge $-1$ anti-instanton 
bundle $\Gamma^\infty(\Sk^7\times_{\SU(2)} \C^2$) on $\S^4$. Similarly 
$\Gamma^\infty(\S^4, \cS^+)$ is isomorphic to the charge $+1$ instanton bundle.

\begin{lma}\label{spbundles} There are the following isomorphisms of right 
$\Cinf(\S^4)$-modules, \begin{align*} \Omega^1(\S^4) &\isom 
\Gamma^\infty(\S^4,\cS^+ \ot \cS^-)\isom 
\Gamma^\infty(\S^4,\cS^+)\ot_{\Cinf(\S^4)} \Gamma^\infty(\S^4,\cS^-) , \\ 
\Omega^0(\S^4) \oplus \Omega^2_-(\S^4) &\isom \Gamma^\infty(\S^4,\cS^- \ot \cS^-) 
\isom \Gamma^\infty(\S^4,\cS^-)\ot_{\Cinf(\S^4)} \Gamma^\infty(\S^4,\cS^-) . 
\end{align*} \end{lma} \begin{proof} Since classically $\Omega^1(S^4) \isom 
\Gamma^\infty(S^4, \cS^+ \otimes \cS^-)$ as $\sigma$-equivariant 
$\Cinf(S^4)$-bimodules, Lemma~\ref{lma:modulesA} shows that $\Omega^1(\S^4) \isom 
\Gamma^\infty(\S^4, \cS^+ \otimes \cS^-)$ as $\Cinf(\S^4)$-bimodules. The 
observations above the Lemma indicate that $\cS^\pm \isom S^7 \times_{\rho^\pm} 
\C^2$ for the spinor representation $\rho^+ \oplus \rho^-$ of $\Spin(4) \isom 
\SU(2) \times \SU(2)$ on $\C^4$, so that \begin{align*} \Gamma^\infty(\S^4, \cS^+ 
\otimes \cS^-) &\isom \Cinf(\Sk^7) \boxtimes_{\rho^+ \otimes \rho^-} (\C^2 
\otimes \C^2)\\ &\isom \left(\Cinf(\Sk^7) \boxtimes_{\rho^+} \C^2\right) 
\otimes_{\Cinf(\S^4)} \left( \Cinf(\Sk^7) \boxtimes_{\rho^-} \C^2 \right) , 
\end{align*} using Proposition~\ref{prop:tensor-modules} in the last line. This 
proves our claim. An analogous statement holds for the second isomorphism. 
\end{proof} Let us forget for the moment the adjoint bundle $\ad(\Sk^7)$. Since 
$\Omega(\S^4) \isom \Omega(S^4)$ as vector spaces and both $\dd$ and the Hodge 
$\ast$ commute with the action of $\bT^2$, the operator $\dd^* + P_- \dd$ can be 
understood as a map from $\Omega^1(S^4) \to \Omega^0(S^4) \op \Omega^2_-(S^4)$ 
(see Section~\ref{section:diff-calc}). Under the isomorphisms of the above Lemma, 
this operator is replaced \cite{AHS78} by a Dirac operator with coefficients in 
$\cS^-$, \be \label{eq:dirac-spinors} D': \Gamma^\infty(\S^4,\cS^+ \ot \cS^-) \to 
\Gamma^\infty(\S^4,\cS^- \ot \cS^-). \ee Twisting by the adjoint bundle 
$\ad(\Sk^7)$, merely results into a composition with the projection $p_{(2)}$ 
defining the bundle $\ad(\Sk^7)$. Hence, eventually the operator $\dd_0^* + 
\dd_1$ is replaced by the Dirac operator \be \label{eq:dirac-twisted} \cD: 
\Gamma^\infty(\S^4, \cS^+ \otimes \cS^-\ot \ad(\Sk^7)) \to 
\Gamma^\infty(\S^4,\cS^- \otimes\cS^- \ot \ad(\Sk^7)), \ee with coefficients in 
the vector bundle $\cS^-\otimes \ad(\Sk^7)$ on $\S^4$.

We have finally arrived to the computation of the index of this Dirac operator 
which we do by means of the Connes-Moscovoci local index formula. It is given by 
the pairing, \be \ind(\cD) = \langle \phi, \chern(\cS^- \ot \ad(\Sk^7)) \rangle = 
\langle \phi, \chern(\cS^-) \cdot \chern(\ad(\Sk^7)) \rangle. \ee In the Appendix 
we recall the expression for both the cyclic cocycle $\phi$ and the Chern 
characters, as well as their realization as operators $\pi_D(\chern (\cdot))$ on 
the the Hilbert space of spinors $\cH$. In \cite{LS04} we computed these 
operators for all modules associated to the principal bundle $\Sk^7 \to \S^4$. In 
particular, for the adjoint bundle we found that \begin{align*} \pi_D 
\big(\chern_0(\ad(\Sk^7))\big) = 3, \qquad \pi_D \big(\chern_1(\ad(\Sk^7))\big) = 
0, \qquad \pi_D \big(\chern_2(\ad(\Sk^7))\big) = 4 (3 \gamma_5). \end{align*} To 
compute the Chern character of the spinor bundle $\cS^-$ we use its mentioned 
identification with the charge $-1$ instanton bundle 
$\Gamma^\infty(\Sk^7\times_{\SU(2)} \C^2$) on $\S^4$. It then follows from 
\cite{CL01} (cf. also \cite{LS04}) that \begin{align*} \pi_D 
\big(\chern_0(\cS^-)\big) = 2, \qquad\pi_D \big(\chern_1(\cS^-)\big) = 0, 
\qquad\pi_D \big(\chern_2(\cS^-)\big) = - 3 \gamma_5. \end{align*} Combining both 
Chern characters and using the local index formula on $\S^4$, we have \be 
\ind(\cD) = 6 ~\resz ~z^{-1} \tr (\gamma_5 |D|^{-2z}) + 0 + \half (2\cdot 4 - 3 
\cdot 1) \resz \tr ( 3 \gamma_5^2 |D|^{-4-2z}), \ee with $D$ identified with the 
classical Dirac operator on $S^4$ (recall that we do not change it in the 
isospectral deformation). Now, the first term vanishes due to $\ind(D)=0$ for 
this classical operator. On the other hand $\gamma_5^2 = \I_4$, and \[ 3 \resz 
\tr ( |D|^{-4-2z}) = 6 \Tr_{\omega}(|D|^{-4}) = 2, \] since the Dixmier trace of 
$|D|^{-m}$ on the $m$-sphere equals $8/m!$ (cf. for instance \cite{GVF01,Lnd97}). 
We conclude that $\ind(\cD) = 5$ and for the moduli space of instantons on 
$\S^4$, we have the following \begin{thm} The tangent space at the base point 
$\nabla_0$ to the moduli space of (irreducible) $\SU(2)$-instantons on $\S^4$ is 
five-dimensional. \end{thm}

\section{Towards Yang--Mills theory on $\M$} \label{section:general}

In this final Section, we shall briefly describe how the Yang--Mills theory on 
$\S^4$ can be generalized to any compact four-dimensional toric noncommutative 
manifold $\M$.

With $G$ a compact semisimple Lie group, let $P \to M$ be a principal $G$ bundle 
on $M$. We take $M$ to be a compact four-dimensional Riemannian manifold equipped 
with an isometrical action $\sigma$ of the torus $\bT^2$. For the construction to 
work, we assume that this action can be lifted to an action $\tilde\sigma$ of a 
cover $\tilde\bT^2$ on $P$ that commutes with the action of $G$. As in 
Section~\ref{section:toric-ncm}, we define the noncommutative algebras 
$\Cinf(\P)$ and $\Cinf(\M)$ as the vector spaces $\Cinf(P)$ and $\Cinf(M)$ with 
star products defined like in \eqref{eq:star-product} with respect to the action 
of $\tilde \bT^2$ and $\bT^2$ respectively or, equivalently as the images of 
$\Cinf(P)$ and $\Cinf(M)$ under the corresponding quantization map $L_\theta$. 
Since the action of $\tilde \bT^2$ is taken to commute with the action of $G$ on 
$P$, the action $\alpha$ of $G$ on the algebra $\Cinf(P)$ given by \be 
\alpha_g(f)(p)=f(g^{-1}\cdot p) \ee induces an action of $G$ by automorphisms on 
the algebra $\Cinf(\P)$. This also means that the inclusion $\Cinf(M) \subset 
\Cinf(P)$ as $G$-invariant elements in $\Cinf(P)$ extends to an inclusion 
$\Cinf(\M) \subset \Cinf(\P)$ of $G$-invariant element in $\Cinf(\P)$. Clearly, 
the action of $G$ translates into a coaction of the Hopf algebra $\Cinf(G)$ on 
$\Cinf(\P)$. \begin{prop} The inclusion $\Cinf(\M) \into \Cinf(\P)$ is a 
(principal) Hopf-Galois $\Cinf(G)$ extension. \end{prop} \begin{proof} One needs 
to establish bijectivity of the canonical map \begin{align*} \chi : \Cinf(\P) 
\ot_{\Cinf(\M)} \Cinf(\P) &\to \Cinf(\P) \otimes \Cinf(G), \\ f' 
\otimes_{\Cinf(\M)} f &\mapsto f' \Delta_R(f) = f' f_{(0)} \otimes f_{(1)}. 
\end{align*} Now, for the undeformed case the bijectivity of the corresponding 
canonical map $\chi^\class:\Cinf(P)\ot_{\Cinf(M)} \Cinf(P) \to \Cinf(P) \otimes 
\Cinf(G)$ follows by the very definition of a principal bundle. Furthermore, 
there is an isomorphism of vector spaces, \begin{align*} T:\Cinf(\P) 
\ot_{\Cinf(\M)} \Cinf(\P) &\to \Cinf(P)\ot_{\Cinf(M)} \Cinf(P) \\ f' 
\ot_{\Cinf(\M)} f &\mapsto \sum_r f'_r \ot_{\Cinf(M)} \tilde\sigma_{r\theta} (f) 
\end{align*} where $f'=\sum_r f'_r$ is the homogeneous decomposition of $f'$ 
under the action of $\tilde \bT^2$. We claim that the canonical map is given as 
the composition $\chi = (L_\theta \otimes \id ) \circ \chi^\class \circ T$; hence 
it is bijective. Indeed, \begin{align*} (L_\theta \otimes \id ) \circ \chi^\class 
\circ T \big( f' \ot_{\Cinf(\M)} f \big) &= \sum_r L_\theta (f'_r 
\tilde\sigma_{r\theta}(f_{(0)})) \ot f_{(1)}\\ &= L_\theta (f' \times_\theta 
f_{(0)}) \ot f_{(1)}=\chi(f' \ot_{\Cinf(\M)} f), \end{align*} since the action of 
$\tilde \bT^2$ on $\Cinf(\P)$ commutes with the coaction of $\Cinf(G)$. 
\end{proof} Noncommutative associated bundles are defined as in \eqref{eq-map} by 
setting $$ \E = \Cinf(\P) \boxtimes_\rho W:= \big\{ f \in \Cinf(\P)\ot W | 
(\alpha_g \ot \id)(f)=(\id \ot \rho(g)^{-1})(f), \, \forall ~g \in G \big\} $$ 
for a representation $\rho$ of $G$ on $W$. These $\Cinf(\M)$ bimodules are finite 
projective since they are of the kind defined in Section~\ref{section:nc-vb} (cf. 
Remark~\ref{rem:associated-modules}). Moreover, Proposition~\ref{prop:end} 
generalizes and reads $\End(\E) \isom \Cinf(\P) \boxtimes_\ad L(W)$, where $\ad$ 
is the adjoint representation of $G$ on $L(W)$. Also, one identifies the adjoint 
bundle as the module coming from the adjoint representation of $G$ on $\g \subset 
L(W)$, namely $\Gamma^\infty(\ad(\P)) := \Cinf(\P)\boxtimes_{ad} \g$.

For a (right) finite projective $\Cinf(\M)$-module $\E$ we define an inner 
product $(\cdot, \cdot)_2$ on $\Hom_{\Cinf(\M)} (\E,\E \ot_{\Cinf(\M)} 
\Omega(\M))$ as in Section~\ref{section:ym}. The Yang--Mills action functional 
for a connection $\nabla$ on $\E$ in terms of its curvature $F$ is then defined 
by \be \YM(\nabla)=(F,F)_2 . \ee This is a gauge invariant, positive and quartic 
functional. The derivation of the Yang--Mills equations \eqref{eq:ym} on $\S^4$ 
does not rely on the specific properties of $\S^4$ and continues to hold on $\M$. 
The same is true for the topological action, and $\YM(\nabla)\geq |\Top(\E)|$ 
with equality iff $\ast_\theta F=\pm F$. In other words, instanton connections 
are minima of the Yang--Mills action.

The explicit construction of instanton connections on $\S^4$ carried over in 
Section~\ref{sect:constrinst} can of course not be generalized to any manifold 
$\M$. Local expressions could in principle be obtained on a ``local chart'' 
$\R^4$ of $\M$ if $\bT^2$ acts on the corresponding local chart $\bR^4$ of $M$. 
On the other hand, the infinitesimal construction of instantons on $\M$ giving 
the dimension of the ``tangent of the moduli space'' can be generalized to any 
toric noncommutative manifold $\M$, again closely following \cite{AHS78}. It is 
essential however the existence of a ``base point'', {i.e.} an instanton 
connection that can be linearly perturbed to obtain a family of infinitesimal 
instantons.

\bigskip \subsection*{Acknowledgements}

We are grateful to an anonymous referee for an excellent review which led to a
much improved version of the paper. We thank Paolo Aschieri and Marc Rieffel for 
useful remarks and suggestions. Part of the work was carried out while GL was 
visiting ESI in Vienna and INI in Cambridge.

\newpage \appendix \section{Local index formula} \label{app:lif}

Suppose in general that $(\A,\cH,D,\gamma)$ is an even $p$-summable spectral 
triple with discrete simple dimension spectrum.  For a projection $e \in M_N(A)$, 
the operator \[ D_e = e(D \ot \II_N)e \] is a Fredholm operator, thought of as 
the Dirac operator with coefficient in the module determined by $e$. The local 
index formula of Connes and Moscovici \cite{CM95} provides a method to compute 
its index via the pairing of suitable cyclic cycles and cocycles. We shall recall 
the ``even'' case since it is the one that is relevant for the present paper.

Let $C_*(\A)$ be the chain complex over the algebra $\A$; in degree $n$, 
$C_n(\A):=\A^{\otimes(n+1)}$. On this complex there are defined the Hochschild 
operator $b:C_n(\A)\to C_{n-1}(\A)$ and the boundary operator $B:C_n(\A)\to 
C_{n+1}(\A)$, satisfying $b^2=0, B^2=0, bB+Bb=0$; thus $(b+B)^2=0$. From general 
homological theory, one defines a bicomplex $CC_*(\A)$ by $CC_{(n,m)}(\A):= 
CC_{n-m}(\A)$ in bi-degree $(n,m)$. Dually, one defines $CC^\ast(\A)$ as 
functionals on $CC_*(\A)$, equipped with the dual Hochschild operator $b$ and 
coboundary operator $B$ (we refer to \cite{C94} and \cite{Lod92} for more details 
on this). \begin{thm}[Connes-Moscovici \cite{CM95}]~\\ \label{thm:cm} \begin{itemize} \item[(a)] 
An even cocycle $\phi^*=\sum_{k\geq0} \phi^{2k}$ in $CC^*(\A)$, $(b+B)\phi^*=0$, is 
defined by the following formul{\ae}. For $k=0$, \be \phi^0(a):=\resz ~z^{-1} 
\tr( \gamma a |D|^{-2z} ); \ee whereas for $k\neq 0$ \be\label{evencocy} 
\phi^{2k}(a^0, \ldots, a^{2k}) := \sum_{\alpha} c_{k,\alpha} \resz ~ \tr \big( 
\gamma a^0 [D,a^1]^{(\alpha_1)} \cdots [D,a^{2k}]^{(\alpha_{2k})} 
|D|^{-2(|\alpha| +k + z)} \big), \ee with \bd c_{k,\alpha}=(-1)^{|\alpha|} 
\Gamma(k+|\alpha|) \big( \alpha! (\alpha_1+1) (\alpha_1+\alpha_2+2) \cdots 
(\alpha_1+ \cdots +\alpha_{2k} + 2k)\big)^{-1} \ed and $T^{(j)}$ denotes the j'th 
iteration of the derivation $T \mapsto [D^2,T]$. \item[(b)] For $e \in K_0(\A)$, 
the Chern character $\chern_*(e) = \sum_{k\geq 0} \chern_k(e)$ is the even cycle 
in $CC_*(\A)$, $(b+B)\chern_*(e)=0$, defined by the following formul{\ae}. For 
$k=0$, \be \chern_0(e):= \tr(e); \ee whereas for $k\neq 0$ \be \label{chernchar} 
\chern_k(e):=(-1)^k \frac{(2k)!}{k!} \sum (e_{i_0 i_1}-\frac{1}{2}\delta_{i_0 
i_1}) \otimes {e_{i_1 i_2} \otimes e_{i_2 i_3} \otimes \cdots \otimes e_{i_{2k} 
i_0} }. \ee \item[(c)] The index of the operator $D_e$ is given by the natural 
pairing between cycles and cocycles, \be \ind D_{e}=\langle 
\phi^\ast,\chern_\ast(e) \rangle . \ee \end{itemize} \end{thm}

For toric noncommutative manifolds, the above local index formula simplifies drastically \cite{LS04}.
\begin{thm}
\label{prop:lif-theta}
For a projection $p \in M_N(\Cinf(\M))$, we 
have $$ \ind D_p = \resz z^{-1} \tr \bigg(\gamma p |D|^{-2z} \bigg)+
\sum_{k \geq 1} \frac{(-1)^k}{k} \resz \tr \bigg( \gamma \big(p-\frac{1}{2} \big) [D,p]^{2k} |D|^{-2(k+z)} 
\bigg)$$ 
and now the trace $\tr$ comprises a matrix trace as well. 
\end{thm}
\begin{proof}
Recall that the quantization map $L_\theta$ preserves the spectral decomposition, for the toric action of $\tilde \bT^n$, of smooth operators (see equation \eqref{spectwist}). 
Then, once extended the deformed $\times_\theta$-product to $\Cinf(\M) \bigcup [D,\Cinf(\M)]$ -- which can be done unambiguously since $D$ is of degree 0 --  we write the local cocycles $\phi^{2k}$ in Theorem~\ref{thm:cm} in terms of the quantization map $L_\theta$:
\begin{multline}
\phi^{2k}\big( L_\theta(f^0),L_\theta(f^1),\ldots, L_\theta(f^{2k}) \big) \\ 
 =  \resz \tr \Big( \gamma L_\theta \big( f^0 \times_\theta [D,f^1]^{(\alpha_1)} \times_\theta \cdots \times_\theta [D,f^{2k}]^{(\alpha_{2k})} \big) 
|D|^{-2(|\alpha| +k + z)}\Big) .
\end{multline}
Suppose now that $f^0,\ldots, f^{2k}\in C^\infty(M)$ are homogeneous of degree $r^0, \ldots, r^{2k}$ respectively, so that the operator 
$f^0 \times_\theta [D,f^1]^{(\alpha_1)} \times_\theta \cdots \times_\theta [D,f^{2k}]$ is a homogeneous element of degree $r=\sum_{j=0}^{2k} r^j$.  
It is in fact a multiple of $f^0 [D,f^1] \cdots  [D,f^{2k}]$ as can be established by working out the $\times_\theta$-products.
Forgetting about this factor -- which is a power of the deformation parameter $\lambda$ -- we obtain from \eqref{twist} that
\be
L_\theta(f^0 \times_\theta [D,f^1]^{(\alpha_1)} \times_\theta \cdots \times_\theta [D,f^{2k}]) = f^0 [D,f^1] \cdots [D,f^{2k}] U(\half r \cdot \theta) .
\ee
Each term in the local index formula for $(C^\infty(\M),\cH,D)$ then takes the form
\bd
\resz  \tr \big( \gamma f^0 [D,f^1]^{(\alpha_1)} \cdots [D,f^{2k}]^{(\alpha_{2k})} |D|^{-2(|\alpha| +k + z)} U(s) \big)
\ed
for $s=\half r \cdot \theta \in \tilde \bT^n$. The appearance of the operator $U(s)$ here is a consequence of the close relation with the index formula for 
$\bT^n$-equivariant Dirac spectral triples. In \cite{CH97} 
an even dimensional compact spin manifold $M$ on which a (connected compact) Lie group $G$ acts by isometries was studied. The equivariant Chern character was defined as an equivariant version of the JLO-cocycle, the latter being an element in equivariant entire cyclic cohomology. The essential point is that  an explicit formula for the above residues was obtained. In the case of a $\bT^n$-action on $M$ this is 
\begin{multline}
\resz  \tr \big( \gamma f^0 [D,f^1]^{(\alpha_1)} \cdots [D,f^{2k}]^{(\alpha_{2k})} |D|^{-2(|\alpha| +k + z)} U(s) \big) \\
=\Gamma(|\alpha|+k) \lim_{t \to 0} t^{|\alpha|+k}  \tr \big( \gamma f^0 [D,f^1]^{(\alpha_1)} \cdots [D,f^{2k}]^{(\alpha_{2k})} e^{-tD^2} U(s) \big) ,
\end{multline}
for every $s \in \tilde \bT^n$. Moreover, from Thm 2 in \cite{CH97} the limit vanishes when $|\alpha|\neq0$. This finishes the proof of our theorem.
\end{proof}

By inserting the symbol $\pi$ for the algebra representation, 
the components of the Chern character are represented as operators on  the Hilbert space $\ch$ by explicit formul{\ae}, 
\be\label{pid}
\pi_D(\chern_k(e)) := (-1)^k \frac{(2k)!}{k!}  \sum (\pi(e_{i_0 i_1})-\frac{1}{2}\delta_{i_0 i_1}) [D, {\pi(e_{i_1 i_2})]
[D, \pi(e_{i_1 i_2})] \cdots [D,\pi(e_{i_{2k} i_0})] } ,
\ee
for $k>0$, while $\pi_D(\chern_0(e)) = \sum \pi(e_{i_0 i_0})$.


\end{document}